%% file: st_article.tex
\documentclass[onefignum,onetabnum]{siamart250211}

\headers{Unfitted Space-Time Discontinuous Galerkin Methods}{F. Heimann, C. Lehrenfeld}

\title{Higher Order Unfitted Discontinuous Galerkin Methods in Space and Time\thanks{Submitted to the editors DATE.
}}

\author{Fabian Heimann\thanks{Department of Mathematics, University College London, UK, \texttt{f.heimann@ucl.ac.uk}} \and
Christoph Lehrenfeld\thanks{Institut für numerische und angewandte Mathematik, Universität Göttingen, Germany, \texttt{lehrenfeld@math.uni-goettingen.de}}}

\input{st_shared_siam}

\usepackage{pifont}
\crefname{assumption}{Assumption}{Assumptions}
\Crefname{assumption}{Assumption}{Assumptions}

\tikzset{external/system call={lualatex \tikzexternalcheckshellescape -extra-mem-top=10000000 -extra-mem-bot=10000000 -halt-on-error -interaction=batchmode -jobname "\image" "\texsource"}}

\ifpdf
\hypersetup{
  pdftitle={Unfitted Space-Time DGDG},
  pdfauthor={F. Heimann, C. Lehrenfeld}
}
\fi

\newcommand{\avg}[1]{\left\{\!\!\left\{ #1 \right\}\!\!\right\}}
\newcommand{\jump}[1]{[ #1 ]}
\newcommand{\fjump}[1]{[\![ #1 ]\!]}

\newcommand{\tripplenorm}[1]{|\!|\!| #1 |\!|\!|}
\newcommand{\tripplejump}[1]{|\!|\![ #1 ]\!|\!|}
\newcommand{\wapp}{\texttt{w}}
\newcommand{\Lapp}{\texttt{L}}
\newcommand{\cdo}{\!\cdot\!}
\newcommand{\n}{\mathbf{n}}
\newcommand{\w}{\mathbf{w}}

\begin{document}

\maketitle

\begin{abstract}
We present and analyse an unfitted space-time discontinuous Galerkin (DG) method of arbitrary order for convection-diffusion problems on moving domains in up to three space dimensions.
The method is discontinuous in space and time and uses upwind fluxes, so that convection is stabilised without additional terms and a Ghost penalty is needed only in the vicinity of the moving boundary.
We prove inf-sup stability in a problem-specific norm, which in addition to the coercive part controls the time derivative in the diffusion-dominated regime and the material derivative in the convection-dominated regime.
From this we derive a priori error estimates of high order in space and time for both regimes, including vanishing diffusion.
Numerical experiments confirm the convergence rates and the robustness with respect to the diffusion coefficient.
\end{abstract}



\section{Introduction}
Partial differential equations (PDEs) on moving domains arise in many applications throughout the natural sciences and engineering, for instance in the transport of dissolved species or surfactants in two-phase flows \cite{grossreusken11, FRACHON2023111734}. Fitted discretisations require the mesh to follow the domain, by mesh deformation or remeshing, which becomes difficult for large deformations or topology changes. \emph{Geometrically unfitted} finite element methods avoid this: the computational mesh is fixed and independent of the geometry, which is described separately, for instance by a level set function \cite{sethian1999level}; see \cite{burman15, Burman_Hansbo_Larson_Zahedi_2025} for overviews.

For unfitted methods on moving domains, the set of active degrees of freedom changes from one time step to the next, so that a method of lines is not directly applicable. Two structurally different approaches have been established. Eulerian time stepping schemes combine finite differences in time with a discrete extension of the solution to a neighbourhood of the current domain \cite{LO_ESAIM_2019}. Space-time formulations discretise the problem variationally on the space-time domain \cite{L_PHD_2015, Z18, BADIA202360, HLP2022}; they need no such extension, allow for arbitrary order in space and time, and remain applicable in the limit of pure transport \cite{burman2025transport}. We follow the second approach.

Within space-time formulations, a further choice is whether the discrete functions are continuous or discontinuous, in time and in space. In time, both options have been studied computationally \cite{HLP2022}, and the numerical analysis has so far been developed for the discontinuous one, combined with continuous functions in space \cite{BADIA202360, heimann2025discretizationerroranalysishigh, burman2025transport}. For this combination, the size of the diffusion coefficient has turned out to matter: the diffusion-dominated regime \cite{heimann2025discretizationerroranalysishigh} and the transport limit \cite{burman2025transport} have been analysed separately, and with different stabilisations. Discontinuous functions in space lead to discontinuous Galerkin (DG) methods, for which upwind fluxes are an established means to treat convection robustly \cite{di2011mathematical, ern2021finite3}.

In this paper we combine the two: the discrete functions are discontinuous in both time and space. The guiding idea is a uniform treatment of the two kinds of discontinuities: a DG discretisation in time amounts to an upwind flux across the time levels, and combining it with upwind fluxes across the spatial facets yields the upwind DG method for the space-time transport operator $\partial_t + \w \cdo \nabla$, supplemented by a symmetric interior penalty term for the diffusion. Convection is thus stabilised by a single mechanism, and the Ghost penalty stabilisation \cite{burman2010ghost}, which is common in unfitted methods, is only needed to control small cuts at the moving boundary. The uniformity carries over to the analysis: inf-sup stability is obtained by one argument for the diffusion-dominated and the convection-dominated regime, in which only the test function differs, the time derivative in the former and the material derivative in the latter. The resulting a priori error bounds are of high order in space and time in both regimes. In the diffusion-dominated regime they are optimal in a norm that contains the diffusion-weighted $H^1$ seminorm, as for the continuous-in-space method \cite{heimann2025discretizationerroranalysishigh}. In the convection-dominated regime they are of order $k + \frac12$ in the $L^2$ norm and of order $k$ for the material derivative, with constants independent of the diffusion coefficient, as is typical for upwind DG methods \cite{di2011mathematical} and as obtained in \cite{burman2025transport} for pure transport.

The analysis assumes that the integrals on the cut space-time domains are evaluated exactly. For the practical realisation by an isoparametric mapping and the analysis of the resulting geometry errors for continuous-in-space methods we refer to \cite{HLP2022, heimann2024geometrically, heimann2025discretizationerroranalysishigh}.

The works closest to the present one stem from the space-time cut finite element method of \cite{HLZ16, Z18}, which is continuous in space and discontinuous in time and has been extended to high order with discrete mass conservation in \cite{MYRBACK2024117245}. A method of the same type, based on cell aggregation instead of a Ghost penalty, is analysed for the diffusion problem in \cite{BADIA202360}, where it is remarked that discontinuous functions in space can be treated as well. Discontinuous functions in both space and time on a cut mesh have been used in \cite{FuFrachonKreissZahedi2022} for linear advection with a moving interface in one space dimension, where conservation and stability of the method are shown.

Both ingredients of our method have counterparts in settings without cut elements in time. On stationary cut meshes, symmetric interior penalty DG methods \cite{GuerkanMassing2019} and upwind DG methods \cite{GuerkanStickoMassing2020} have been analysed with Ghost penalty stabilisation; for hyperbolic problems with explicit time stepping, small cut cells are treated by Ghost penalties on macro-elements \cite{FuKreissZahedi2024, FuKreissXinZahedi2026} or by the domain-of-dependence stabilisation \cite{EngwerMayNuessingStreitbuerger2020}. On fitted meshes that move with the domain, space-time DG methods that are discontinuous in space and time have been analysed for advection-diffusion problems in \cite{SudirhamVegtDamme2006} and, in hybridised form, in \cite{KirkHorvathCesmeliogluRhebergen2019, WangRhebergen2025}, the latter for the advection-dominated regime. In comparison to these works, we consider convection-diffusion problems on moving domains with an unfitted mesh in up to three space dimensions and provide an a priori error analysis in both regimes.

Beyond the scalar model problem considered here, the fully discontinuous setting is of interest as a starting point for problems for which DG discretisations are a standard choice, such as hyperbolic systems and nonlinear conservation laws on moving domains, where local conservation and upwind-type numerical fluxes play a central role. These extensions are left for future work.

\subsection*{Structure of the paper}
The remainder of the paper is structured as follows. In \cref{sec:model_problem} we introduce the model problem. \Cref{sec:method} introduces the unfitted space-time DG method including the Ghost penalty stabilisation. The stability analysis, which distinguishes a diffusion dominated and a convection dominated regime, is carried out in \cref{sect_stab_analysis}, and the a priori error analysis in \cref{sect_a_pr_est}. Numerical examples are presented in \cref{sec:numex}.

\section{Model problem} \label{sec:model_problem}
We consider the convection-diffusion equation on a moving domain $\Omega(t) \subset \mathbb{R}^d$, $d \in \{1,2,3\}$, which is smooth and bounded for all $t \in [0,T]$, $T > 0$: Find $u = u(x,t)$ such that
\begin{subequations}
 \label{strongformproblem}
\begin{align}
 \partial_t u + \w \cdot \nabla u - \alpha \Delta u &= f \quad \text{in } \Omega(t), \ t \in (0,T], \label{strongformproblem:1} \\
 \nabla u \cdot \n_{\partial \Omega} &= 0 \quad \text{on } \partial \Omega(t), \ t \in (0,T], \label{strongformproblem:2} ~ (\textnormal {if } \alpha > 0) ,\\
 u(\cdot,0) &= u_0 \quad \text{in } \Omega(0). \label{strongformproblem:3}
\end{align}
\end{subequations}
Here, $u$ is an unknown scalar field, e.g. a chemical concentration, $\alpha \geq 0$ is the constant diffusion coefficient (which is allowed to vanish), $f$ is a source term, $u_0$ are the initial data and $\w = \w(x,t)$ is a given convection velocity. We assume that $\w$ is given on a fixed background domain $\tilde \Omega \supset \Omega(t)$, $t \in [0,T]$, with $\w \in W^{1,\infty}(\tilde \Omega \times (0,T))$, is \emph{divergence-free}, $\operatorname{div} \w = 0$, and governs the motion of the domain: the normal velocity $V_{\partial \Omega}$ of the boundary satisfies
\begin{equation}
 V_{\partial \Omega} = \w \cdo \n_{\partial \Omega} \quad \text{on } \partial \Omega(t), \ t \in (0,T], \label{eq:domain_motion}
\end{equation}
where $\n_{\partial \Omega}$ denotes the outer unit normal of $\Omega(t)$. Together with \cref{eq:domain_motion}, the boundary condition \cref{strongformproblem:2} ensures that there is no flux across the moving boundary. In the case of vanishing diffusion, $\alpha = 0$, the problem is complete without a boundary condition, as the domain moves with the same velocity $\w$ as the concentration field.

\section{Unfitted space-time DG formulation} \label{sec:method}

\subsection{Preliminaries and notation}
\paragraph{Domains and meshes}
With the background domain $\tilde \Omega \subset \mathbb{R}^d$ from \cref{sec:model_problem}, which contains $\Omega(t)$ for all $t \in [0,T]$, let $\phi$ be a level set function describing the moving domain, $\Omega(t) = \{ x \in \tilde \Omega \mid \phi(x,t) < 0 \}$. The space-time domain is
\begin{equation}
 Q := \bigcup_{t \in (0,T)} \Omega(t) \times \{t\} = \{ (x,t) \in \tilde \Omega \times (0,T) \, | \, \phi(x,t) < 0 \} \subset \tilde \Omega \times (0,T).
\end{equation}

Let ${\tilde{\mathcal{T}}}_h  = \{T\}$ be a shape-regular and quasi-uniform simplicial mesh of $\tilde \Omega$, i.e. consisting of intervals, triangles or tetrahedra for $d = 1, 2, 3$, with mesh size $h = \max_{T \in {\tilde{\mathcal{T}}}_h} h_T$, $h_T = \operatorname{diam}(T)$, and set of facets ${\tilde{\mathcal{F}}}_h$.
The time interval $[0,T]$ is divided into uniform steps of size $\Delta t$ with $I_n = (t_{n-1},t_n]$, $n=1,\dots,N$, $t_0 = 0$, $t_N = T$.
The tensor product of ${\tilde{\mathcal{T}}}_h$ with a time interval yields the space-time mesh ${\tilde{ \mathcal{Q}}}_T^n = \{Q_T^n \}_{T\in {\tilde{ \mathcal{T}}}_h}$ of the time slab $\tilde Q^n := \tilde \Omega \times I_n$, consisting of the space-time prisms $Q_T^n = T \times I_n$, and ${\tilde{ \mathcal{Q}}}_T = \bigcup_{n=1}^N {\tilde{ \mathcal{Q}}}_T^n$.

In the discretisation only the space-time elements that contribute to the space-time domain $Q$ are used. With the physical part of each time slab, $Q^n := Q \cap \tilde Q^n$, we define the \emph{active} space-time triangulation and the corresponding (staircase) domain slab by slab as
\begin{equation}
 \mathcal{Q}_T^n := \{ Q_T^n \in {\tilde{\mathcal{Q}}}_T^n \, | \, Q_T^n \cap Q^n \neq \emptyset \} \subset {\tilde{\mathcal{Q}}}_T^n, \qquad Q_{\mathcal{E}}^n := \bigcup_{Q_T^n \in \mathcal{Q}_T^n} Q_T^n \supset Q^n, \label{eq:def_active_mesh}
\end{equation}
and globally as $\mathcal{Q}_T := \bigcup_{n=1}^N \mathcal{Q}_T^n$ and $Q_{\mathcal{E}} := \bigcup_{n=1}^N Q_{\mathcal{E}}^n$.

On each time slab we define the set of active facets as those that intersect $Q^n$ and the corresponding space-time facets,
$$
\mathcal{F}_h^n = \{ F \in {\tilde{\mathcal{F}}}_h \mid (F \times I_n) \cap Q^n \neq \emptyset \}, \qquad \mathcal{Q}_F^n = \{ F \times I_n \mid F \in \mathcal{F}_h^n  \},
$$
and $\mathcal{F}_h$ and $\mathcal{Q}_F$ as the respective unions over all time slabs.
Finally, we set $\Omega^n := \Omega(t_n)$, $n = 0, \dots, N$, for the physical domain at the time levels.

\paragraph{Discrete spaces}
With the spatial and temporal polynomial orders $k_s$ and $k_t$ we define the spatial DG space and the fully discontinuous space-time spaces on the background mesh,
\begin{align*}
 \tilde V_h^s &:= \{ v \in L^2(\tilde \Omega) \, | \, v|_T \in \mathcal{P}^{k_s}(T) \ \forall T \in \tilde{\mathcal{T}}_h \}, \qquad \tilde W_h^n := \tilde V_h^s \otimes \mathcal{P}^{k_t}(I_n), \\
 \tilde W_h &:= \{ v \in L^2(\tilde \Omega \times (0,T)) \mid v|_{\tilde Q^n} \in \tilde W_h^n, \ n = 1, \dots, N \}.
\end{align*}
The discrete spaces of the method are their restrictions to the active elements,
\begin{equation}
 W_h^n := \{ v|_{Q_{\mathcal{E}}^n} \mid v \in \tilde W_h^n \}, \qquad W_h := \{ v|_{Q_{\mathcal{E}}} \mid v \in \tilde W_h \}.
\end{equation}
Note that discrete functions in $W_h$ are allowed to be discontinuous across the time slab boundaries.

\paragraph{Inner products}
With $(\cdot,\cdot)_S$ we denote the $L^2$ inner product on a domain $S$ and with $\Vert \cdot \Vert_S$ the corresponding norm.
For sets of elements and facets, inner products are understood as sums. On each time slab we define
$(u,v)_{\mathcal{Q}_T^n} := \sum_{Q_T^n \in \mathcal{Q}_T^n} (u,v)_{Q_T^n}$
and
$(u,v)_{Q^n} := \sum_{Q_T^n \in \mathcal{Q}_T^n} (u,v)_{Q_T^n \cap Q^n}$,
and for the space-time facets \emph{within} each time slab
$(u,v)_{\mathcal{Q}_F^n} := \sum_{F \in \mathcal{F}_h^n} (u,v)_{F \times I_n}$ and
$(u,v)_{\mathcal{Q}_F^n \cap Q} := \sum_{F \in \mathcal{F}_h^n} (u,v)_{(F \times I_n) \cap Q}$.
Omitting the index $n$ denotes the sum over all time slabs, e.g. $(u,v)_Q := \sum_{n=1}^N (u,v)_{Q^n}$ or $(u,v)_{\mathcal{Q}_F} := \sum_{n=1}^N (u,v)_{\mathcal{Q}_F^n}$.
For the union of all space-time facets, i.e. the spatial facets within the time slabs and the time slab boundaries, we use the notation $(u,v)_{\mathcal{Q}_F^{\text{st}}} := (u,v)_{\mathcal{Q}_F} + \sum_{n=0}^{N-1} (u,v)_{\Omega^n}$ and $(u,v)_{\mathcal{Q}_F^{\text{st}} \cap Q} := (u,v)_{\mathcal{Q}_F \cap Q} + \sum_{n=0}^{N-1} (u,v)_{\Omega^n}$ (the time slab boundaries $\Omega^n$ are physical domains already).
The final time level $\Omega^N$, on which no (upwind) jump is defined, is excluded here and will always be treated separately. With $\partial_s$ we denote the \emph{spatial} part of a boundary.
All bilinear forms of the discretisation below are formulated on the physical domain and the physical parts of the facets, i.e. on $Q^n$ and $\mathcal{Q}_F^n \cap Q$, whereas the full elements $\mathcal{Q}_T^n$ and facets $\mathcal{Q}_F^n$ only enter through the Ghost penalty stabilisation and in the analysis.

\paragraph{Space-time normals, jumps and averages}
We define $\w^\ast = (\w,1)$ as the space-time velocity and $\nabla^\ast = (\nabla,\partial_t)$ as the space-time gradient so that the material derivative also reads as
$\partial_t  + \w \cdo \nabla= \w^\ast \cdo \nabla^\ast$.

On and across space-time facets, i.e. spatial facets and time slab boundaries we will introduce further notation for the upcoming DG operators.
On each space-time facet we define the unit space-time normal $\n^\ast$ so that $\w^\ast \cdot \n^\ast \geq 0$. The spatial facet normals will accordingly be defined so that $\n^\ast = (\n,0)$. Note that this means that the normal always points in the direction of the spatial velocity on spatial facets and in positive time direction on time slab boundaries. We then define the space-time downwind and upwind values, their jump and average as
$$
 u_{\pm} = \lim_{s \searrow 0} u((x,t)  \pm s \n^\ast), \qquad \fjump{u} = u_+ - u_-,
 \qquad \avg{u} = \frac{u_+ + u_-}{2}.
$$
Note that the orientation of $\n^\ast$ is defined pointwise, so that on a spatial facet $F$ the upwind side may change in space and time. To be precise, fix an orientation $\n_F$ of $F$, pointing from $T_1$ to $T_2$, and set $s_F := 1$ where $\w \cdo \n_F \geq 0$ and $s_F := -1$ where $\w \cdo \n_F < 0$; then $\n := s_F \n_F$ (so that where the normal velocity vanishes, the fixed orientation is used) and, with the fixed-orientation jump $\delta_F u := u|_{T_2} - u|_{T_1}$, we have $\fjump{u} = s_F \, \delta_F u$, in particular $|\fjump{u}| = |\delta_F u|$ pointwise. The latter will be used in the analysis whenever an operator acting in time (and hence not commuting with the pointwise orientation) is applied to a jump.

\subsection{Discrete formulation}
We use a symmetric interior penalty discretisation for the (spatial) diffusion and upwind fluxes for the space-time convection, i.e. for the convective and the time derivative alike. We first formulate the problem with global bilinear and linear forms $B_h(u,v)$, $f_h(v)$, which is the form used in the analysis, and extract the corresponding time stepping formulation in \cref{rem:timestepping} below.

\subsubsection*{Global discrete formulation}
We divide the bilinear form $B_h$ into several contributions, each of which is the sum of its slab-wise counterparts, $B_h := C + A + U_h + N_h$,
\begin{align*}
  C := &\sum_{n=1}^N C^n, \quad A := \sum_{n=1}^N A^n, \quad
  U_h := \sum_{n=1}^N U^n_h, \quad N_h := \sum_{n=1}^N N^n_h, \quad \textnormal{ with } \\
 C^n(u,v) =&~  (\partial_t u + \w \cdot \nabla u, v)_{Q^n}, ~~ A^n(u,v) = \alpha (\nabla u, \nabla v)_{Q^n}, \\
 U_h^n(u,v) =&~ (|\w \cdot \n| \fjump{u},{v_+})_{\mathcal{Q}_F^n \cap Q} + (\fjump{u}, v_+)_{\Omega^{n-1}}, \\
 N^n_h(u,v) =&~ (\avg{\alpha \nabla u \cdo \n}\!, \fjump{v})_{\mathcal{Q}_F^n \cap Q} + (\avg{\alpha \nabla v \cdo \n}\!,\fjump{u})_{\mathcal{Q}_F^n \cap Q} + \frac{\alpha \lambda}{h} (\fjump{u}, \fjump{v})_{\mathcal{Q}_F^n \cap Q},
\end{align*}
where we set $\fjump{u} = u_+$ (i.e. $u_- := 0$) on $\Omega^0$ and denote $\lambda > 0$ as the interior penalty parameter. Note that with the convention $\fjump{u} = u_+ - u_-$, $\n$ pointing from the upwind to the downwind side, the consistency terms in $N_h^n$ carry a positive sign, cf. the proof of \cref{lemma_consistency} below.
Here, $C$ and $U_h$ are the volume and facet parts of the space-time convection (the material derivative), and $A$ and $N_h$ those of the spatial diffusion. The Ghost penalty stabilisation bilinear form $J(\cdot,\cdot)$ is defined in \cref{eq:def_J} below with the same slab-wise structure.

\begin{remark}[Space-time upwinding] \label[remark]{rem:st_upwind}
  Note that we can also rewrite $C$ and $U_h$ in a more direct way in the sense of space-time operators:
  $$
  (C + U_h)(u,v)
  =
  (\w^\ast \cdo \nabla^\ast u, v)_{Q}
  +
  ( |\w^\ast \cdo \n^\ast| \fjump{u}, v_+)_{\mathcal{Q}_F^{\text{st}} \cap Q}.
  $$
  This illustrates that the upwind mechanism to handle the convection in space is the same as the transfer of information in the time stepping scheme.
\end{remark}

The right-hand side $f_h$ is given as
 \begin{equation}
  f_h(v) := (f,v)_{Q} + (u_0, v_+)_{\Omega^0}.
 \end{equation}

 The discretisation is complemented by a Ghost penalty stabilisation bilinear form. We use a space-time version of the direct Ghost penalty as in \cite{preuss18, heimann20, HLP2022, heimann2025discretizationerroranalysishigh, burman2025transport}.
 It operates on interior facets of the background mesh. Fix any such $F \in \tilde{\mathcal{F}}_h$ and let $T_1, T_2 \in \tilde{\mathcal{T}}_h$ be the two elements such that $F = \overline{T_1} \cap \overline{T_2}$. We denote the union of those elements as the facet patch $\omega_F = T_1 \cup T_2$. On this facet patch, we define a (volumetric) jump operator for elementwise polynomials as
\begin{equation}
 \jump{u}_{\omega_F} |_{T_i} (x,t) = u|_{T_i}(x,t) - \mathcal{E}^p (u|_{T_{3-i}})(x,t), \quad i=1,2, \label{eq:def_gp_jump}
\end{equation}
where $\mathcal{E}^p \colon \mathcal{P}^{k_s}(T_{3-i}) \otimes \mathcal{P}^{k_t}(I_n) \to \mathcal{P}^{k_s}(\mathbb{R}^d) \otimes \mathcal{P}^{k_t}(I_n)$ denotes the canonical polynomial extension (in space). In the pure transport case \cite{burman2025transport} the Ghost penalty acts on all facets between active elements, as it also serves as the stabilisation of the convection (in the spirit of continuous interior penalty methods). Here, the upwind fluxes take over this role and the Ghost penalty is only needed to handle the cut elements. Accordingly, the stabilisation operates only on the facets of the cut region within one time slab, i.e. on facets between active elements of which at least one is not contained in the physical space-time domain,
\begin{align}
  \mathcal{F}^{n}_{\text{gp}} := \{F = \overline{T_1} \! \cap\! \overline{T_2} \in \tilde{\mathcal{F}}_h \textnormal{ s.t. } T_1 \neq T_2 \in \tilde{\mathcal{T}}_h, \ Q_{T_1}^n, Q_{T_2}^n \in \mathcal{Q}_T^n, \ Q_{T_1}^n \cup Q_{T_2}^n \not\subseteq \overline{Q^n} \}. \label{eq:def_gp_facets}
 \end{align}
Note that an element with $Q_T^n \not\subseteq \overline{Q^n}$ is not necessarily cut at a fixed time, but may also be an element that is swept by the moving boundary during the time slab. We assume throughout that $Q$ has a positive distance to $\partial \tilde \Omega$, so that all active facets are interior facets, and we split the active facets into the stabilised ones and the remaining \emph{interior} facets,
\begin{equation}
 \mathcal{F}_h^n = (\mathcal{F}_h^n \cap \mathcal{F}^n_{\text{gp}}) \, \dot\cup \, \mathcal{F}^n_{\text{int}}, \qquad \mathcal{F}^n_{\text{int}} := \mathcal{F}_h^n \setminus \mathcal{F}^n_{\text{gp}}, \qquad F \times I_n \subseteq \overline{Q^n} \quad \text{for all } F \in \mathcal{F}^n_{\text{int}}, \label{eq:def_int_facets}
\end{equation}
where the last property holds as both elements adjacent to $F \in \mathcal{F}^n_{\text{int}}$ are contained in $\overline{Q^n}$. In other words, every active facet that is cut by the lateral boundary of $Q$ within the time slab is stabilised, whereas the interior facets are not. Correspondingly, we write $\mathcal{Q}_F^{\text{gp}}$ and $\mathcal{Q}_F^{\text{int}}$ for the union over all time slabs of the space-time facets $F \times I_n$ with $F \in \mathcal{F}_h^n \cap \mathcal{F}^n_{\text{gp}}$ and $F \in \mathcal{F}^n_{\text{int}}$, respectively, so that $\mathcal{Q}_F = \mathcal{Q}_F^{\text{gp}} \cup \mathcal{Q}_F^{\text{int}}$ and $\mathcal{Q}_F^{\text{int}} \subseteq \overline{Q}$.
Based on this notation, we define the stabilisation bilinear forms
 \begin{align}
 j_F^n(u,v) &:= \int_{t_{n-1}}^{t_n} \int_{\omega_F} \jump{u}_{\omega_F} \jump{v}_{\omega_F} \, \mathrm{d}x \, \mathrm{d}t,
 \qquad
 j^n_h(u,v) := \sum_{F \in \mathcal{F}^{n}_{\text{gp}}} j_F^n(u,v), \label{eq:def_jhn} \\
 J(u,v) &:= c_J \sum_{n=1}^N j^n_h(u,v), \qquad c_J := \gamma_J \left( \frac{\alpha}{h^2} + \frac{\Vert \w \Vert_\infty}{h} \right), \label{eq:def_J}
\end{align}
where $\gamma_J > 0$ denotes a scalar valued stabilisation parameter. The scaling $c_J$ combines the scaling $h^{-2}$ that is typical for the diffusion operator, cf. \cite{preuss18, heimann2025discretizationerroranalysishigh}, with the scaling $h^{-1}$ that has been used for the pure transport problem in \cite{burman2025transport}. Note that $c_J \simeq \alpha h^{-2} + \Vert \w \Vert_\infty h^{-1}$ is -- up to the constant $\gamma_J$ -- the natural scaling of the bilinear form $B_h$ itself w.r.t. $h$. Here and in the following $\Vert \w \Vert_\infty := \Vert \w \Vert_{L^\infty(\tilde \Omega \times (0,T))}$.

We can now pose the discrete problem: Find $u_h \in W_h$ such that for all $v_h \in W_h$
\begin{equation}
 B_h(u_h,v_h) + J(u_h,v_h) = f_h(v_h). \label{eq_discrete_problem}
\end{equation}

\begin{remark}[Mass conserving formulation] \label[remark]{rem:Bmc}
 By partial integration (in space and time) one can also formulate the convective part of $B_h$ in a \emph{mass conserving} form: testing with $v = 1$ on one time slab yields $(u_-,1)_{\Omega^n} - (u_-,1)_{\Omega^{n-1}} = (f,1)_{Q^n}$ for the discrete solution (with $u_-$ read as $u_0$ on $\Omega^0$).
 We will make use of this form in the analysis, see \cref{lemma_B_Bmc}, as it is the key for the symmetric testing. As we consider exact geometry handling in this paper, the two forms coincide. For discrete geometries this is no longer the case, see \cite{heimann2025discretizationerroranalysishigh, burman2025transport} for the resulting (consistency) errors.
\end{remark}

\begin{remark}[Time stepping formulation] \label[remark]{rem:timestepping}
 As discrete functions are discontinuous across the time slab boundaries and all bilinear forms are sums of slab-wise contributions, the only coupling between the time slabs in \cref{eq_discrete_problem} is through the temporal upwind term $(\fjump{u}, v_+)_{\Omega^{n-1}} = (u_+, v_+)_{\Omega^{n-1}} - (u_-, v_+)_{\Omega^{n-1}}$ in $U_h^n$, in which $u_-$ is the trace of the solution on the previous time slab. Testing \cref{eq_discrete_problem} with functions $v_h \in W_h^n$ (extended by zero to the other time slabs) hence shows that it is equivalent to the following time stepping scheme: Set $u^0_- := u_0$ and for $n = 1, \dots, N$ find $u^n \in W_h^n$ such that for all $v \in W_h^n$
 \begin{equation}
  B_h^n(u^n, v) + c_J \, j_h^n(u^n, v) = (f, v)_{Q^n} + (u^{n-1}_-, v_+)_{\Omega^{n-1}}, \label{eq_time_stepping}
 \end{equation}
 where $B_h^n := C^n + A^n + U_h^n + N_h^n$ with the convention $\fjump{u} = u_+$ on $\Omega^{n-1}$ in $U_h^n$ (as for $n=1$ above), i.e. the upwind value from the previous time slab is moved to the right hand side, and $u^{n-1}_-$ denotes the trace of $u^{n-1}$ on $\Omega^{n-1}$ for $n \geq 2$. The solution of \cref{eq_discrete_problem} is then given by $u_h|_{\tilde Q^n} = u^n$. Note that the active mesh changes from one time slab to the next; $u^{n-1}_-$ is only needed on $\Omega^{n-1}$, which is covered by both $\mathcal{Q}_T^{n-1}$ and $\mathcal{Q}_T^n$.
\end{remark}
\section{Discrete stability and continuity analysis} \label{sect_stab_analysis}

As part of the a priori error analysis we treat the stability in this section. The stability analysis splits into two steps: First, in a coercivity-based approach we obtain control on a (semi-)norm by symmetric testing, i.e. by testing with (a temporally weighted version of) $u_h$ itself. However, especially in the case of vanishing diffusion this semi-norm is too weak to obtain reasonable error bounds. Hence, in a second step we consider testing with a properly scaled derivative as usual in the analysis of upwind DG methods, cf. \cite{preuss18, heimann2025discretizationerroranalysishigh, burman2025transport} for the space-time context. This yields control in a stronger, but discrete, norm. Combining both test functions yields inf-sup stability.

\subsection{Preliminaries and assumptions} \label{ssec:prelim}
Throughout, we write $\lesssim$ ($\gtrsim$) to indicate inequalities up to a constant independent of $h$, $\Delta t$, $\alpha$, $\w$, or the location of the cut within elements, and $a \simeq b$ if $a \lesssim b$ and $b \lesssim a$. Constants may depend on the shape regularity of $\tilde{\mathcal{T}}_h$, the polynomial degrees $k_s, k_t$, the final time $T$, the weighting parameter $\sigma$ introduced in \cref{ssec:coercivity}, the interior penalty parameter $\lambda$, the Ghost penalty parameter $\gamma_J$ and the constants in \cref{ass:gp} below. Moreover, we assume $\Delta t \leq 1$ throughout, which allows us to bound powers of $\Delta t$ by lower powers; this is no restriction if $T \leq 1$ and a mild one otherwise. Likewise, we assume $h \lesssim 1$ throughout (which holds as $\tilde \Omega$ is bounded), so that powers of $h$ can be bounded by lower powers.
To keep the presentation readable, we tacitly assume $\Vert \w \Vert_\infty \simeq 1$ in the sense that we do not track the dependence of constants on $\Vert \w \Vert_{W^{1,\infty}}$; we do, however, keep track of the dependence on $\alpha$, as the regime $\alpha \to 0$ is of particular interest.

We will frequently make use of the following standard discrete inequalities for $v \in W_h$, $Q_T^n = T \times I_n \in \mathcal{Q}_T^n$, $F \subset \partial T$:
\begin{subequations} \label{eq:std_discrete_ineqs}
\begin{align}
 \Vert v \Vert_{F \times I_n}^2 &\lesssim h^{-1} \Vert v \Vert_{Q_T^n}^2, \qquad
 \Vert v_\pm \Vert_{T}^2 \lesssim \Delta t^{-1} \Vert v \Vert_{Q_T^n}^2, \label{eq:std_trace} \\
 \Vert \nabla v \Vert_{Q_T^n} &\lesssim h^{-1} \Vert v \Vert_{Q_T^n}, \qquad
 \Vert \partial_t v \Vert_{Q_T^n} \lesssim \Delta t^{-1} \Vert v \Vert_{Q_T^n}, \label{eq:std_inverse}
\end{align}
and for $v \in H^1(Q_T^n)$, $v_\pm$ denoting the traces at $t_{n-1}$ and $t_n$, respectively,
\begin{align}
 \Vert v \Vert_{F \times I_n}^2 &\lesssim h^{-1} \Vert v \Vert_{Q_T^n}^2 + h \Vert \nabla v \Vert_{Q_T^n}^2, \qquad
 \Vert v_\pm \Vert_{T}^2 \lesssim \Delta t^{-1} \Vert v \Vert_{Q_T^n}^2 + \Delta t \Vert \partial_t v \Vert_{Q_T^n}^2. \label{eq:cont_trace}
\end{align}
\end{subequations}
Note that all these inequalities are formulated on the full (uncut) space-time prisms $Q_T^n$ and facets $F \times I_n$. The role of the Ghost penalty stabilisation is to transfer control from the physical domain $Q^n$ to the full active mesh region $\mathcal{Q}_T^n$ so that the inequalities become applicable.  Moreover, the spatial inequalities in \cref{eq:std_trace} and \cref{eq:std_inverse} only use that $v(\cdot,t)$ is a polynomial on $T$ for each $t \in I_n$; they hence also hold for functions that are piecewise polynomial in space without being polynomial in time, such as $\varphi v$ and $\Pi_t(\varphi v)$ for $v \in W_h$ and $\varphi = \varphi(t)$.

\begin{figure}
 \tikzset{external/export next=false}
 \begin{minipage}[c]{0.55\textwidth}
 \begin{tikzpicture}
  \begin{loglogaxis}[
    scale only axis, width=0.78\linewidth, height=0.6\linewidth,
    xmin=1e-4, xmax=1, ymin=1e-5, ymax=1,
    xtick={1e-4,1e-3,1e-2,1e-1,1}, ytick={1e-5,1e-4,1e-3,1e-2,1e-1,1},
    axis on top, clip=true, tick label style={font=\scriptsize}, label style={font=\small},
   ]
   \addplot[fill=blue!25, draw=blue!70!black, thick] coordinates {(1e-4,1e-2) (1e-2,1e-2) (3.1623e-4,1e-5) (1e-4,1e-5)} --cycle;
   \addplot[fill=white, draw=none] coordinates {(1e-4,1e-3) (1e-3,1e-2) (1e-4,1e-2)} --cycle;
   \addplot[pattern=crosshatch, pattern color=red!60!black, draw=none] coordinates {(1e-4,1e-3) (1e-3,1e-2) (1e-4,1e-2)} --cycle;
   \addplot[fill=orange!45, draw=orange!80!black, thick, opacity=0.8] coordinates {(1e-4,0.5e-4) (1,0.5) (1,1) (0.5,1) (1e-4,2e-4)} --cycle;
   \addplot[black, dashed, thick] coordinates {(3.1623e-4,1e-5) (1e-2,1e-2) (1,1)};
   \draw[red!70!black, thick, dash dot] (axis cs:1e-4,1e-3) -- node[pos=0.8, sloped, below, font=\scriptsize, inner sep=2pt] {$v_\partial \Delta t = h$} (axis cs:1e-1,1);
   \draw[-latex, red!70!black, thick] (axis cs:3e-3,3e-2) -- node[pos=1, left, font=\scriptsize, inner sep=1pt, text=red!70!black] {$v_\partial \to 0$} (axis cs:1.5e-3,6e-2);
   \addplot[only marks, mark=*, mark size=1.6pt] coordinates {(1e-2,1e-2)};
   \node[font=\scriptsize, anchor=north west] at (axis cs:1.1e-2,1e-2) {$(\alpha,\alpha)$};
   \draw[-latex, thick] (axis cs:6e-3,6e-3) -- node[pos=0.5, below right, font=\scriptsize, inner sep=1pt, fill=white, fill opacity=0.85, text opacity=1] {$\alpha \to 0$} (axis cs:1.2e-3,1.2e-3);
   \node[font=\small] at (axis cs:1.8e-4,2.5e-5) {D};
   \node[font=\small] at (axis cs:7e-4,2.5e-3) {D};
   \node[font=\small] at (axis cs:3e-1,5e-2) {C};
   \node[font=\scriptsize, align=center] at (axis cs:6e-4,3e-1) {not covered\\($\Delta t > \alpha$, $\Delta t \gg h$)};
   \node[font=\scriptsize, align=center] at (axis cs:5e-2,1.5e-4) {not covered\\($\Delta t \ll \min\{h, h^2/\alpha\}$)};
  \end{loglogaxis}
  \node[anchor=west, font=\small, inner sep=1pt] at ([xshift=10pt, yshift=-8pt]current axis.south east) {$h$};
  \node[anchor=south, font=\small, inner sep=1pt] at ([xshift=-11pt, yshift=6pt]current axis.north west) {$\Delta t$};
 \end{tikzpicture}
 \end{minipage}\hfill
 \begin{minipage}[c]{0.42\textwidth}
 \caption{Admissible $(h,\Delta t)$ regimes for $\Vert \w \Vert_\infty \simeq 1$, $\alpha = 10^{-2}$ (constants set to one). Blue: Case D \cref{eq:ass_caseD}, $h^2/\alpha \lesssim \Delta t \lesssim \alpha$. Orange: Case C \cref{eq:ass_caseC}, $\Delta t \simeq h$, any $\alpha \geq 0$. Dashed: lower bound on $\Delta t$ from \cref{eq:ass_gp_scaling}. Dash-dotted: boundary CFL condition of \cref{ass:gp} for $v_\partial = 10^{-1}$, which excludes the cross-hatched part.}
 \label{fig:regimes}
 \end{minipage}
\end{figure}
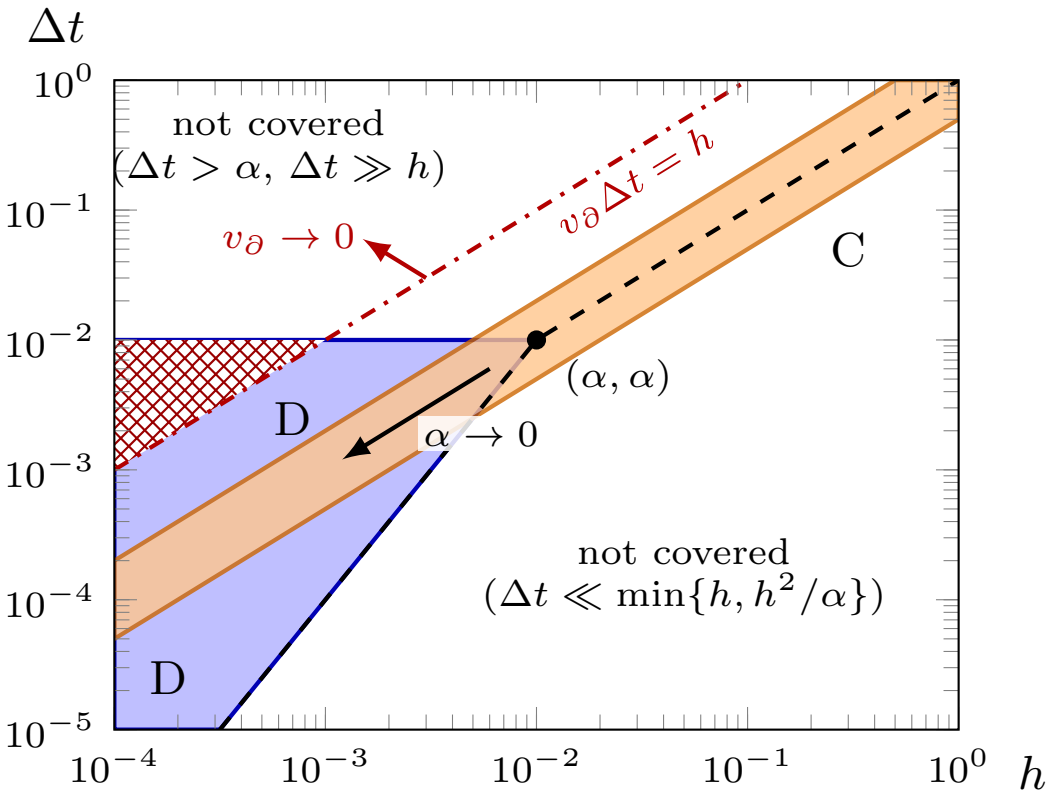

\paragraph{Two regimes}
In the stability analysis we distinguish two regimes, which differ in the derivative that is controlled in the discrete norm and in the assumptions on the relation between $h$, $\Delta t$ and the problem parameters $\alpha$ and $\w$:
\begin{description}
 \item[Case D (diffusion dominated)] We choose $\mathcal{L} u := \partial_t u$ and assume
 \begin{equation}
  h^2 \lesssim \alpha \Delta t \quad \text{ and } \quad \Vert \w \Vert_\infty^2 \Delta t \lesssim \alpha. \label{eq:ass_caseD}
 \end{equation}
 \item[Case C (convection dominated)] We choose $\mathcal{L} u := \partial_t u + \w \cdo \nabla u = \w^\ast \cdo \nabla^\ast u$ and assume
 \begin{equation}
  h \simeq \Vert \w \Vert_\infty \Delta t, \quad \text{i.e. } \Vert \w \Vert_\infty \Delta t \lesssim h \text{ (CFL-type condition) and } h \lesssim \Vert \w \Vert_\infty \Delta t. \label{eq:ass_caseC}
 \end{equation}
 No assumption on $\alpha \geq 0$ is made; in particular $\alpha = 0$ is allowed.
\end{description}
The first assumption in \cref{eq:ass_caseD} is the (inverse CFL-type) assumption that has also been used in \cite{heimann2025discretizationerroranalysishigh} for $\alpha=1$; the second one states that the convection does not dominate the problem on the scale of one time step. The assumption \cref{eq:ass_caseC} is the one of \cite{burman2025transport} and expresses that space and time are refined simultaneously. In both cases, the Ghost penalty scaling satisfies
\begin{equation}
 \Delta t^{-1} \lesssim \alpha h^{-2} + \Vert \w \Vert_\infty h^{-1} \simeq c_J, \label{eq:ass_gp_scaling}
\end{equation}
which is the property that is needed to control temporal derivatives in the Ghost penalty norm. Moreover, in Case D the assumptions \cref{eq:ass_caseD} imply $\Vert \w \Vert_\infty h \lesssim \alpha$, so that $c_J \simeq \alpha h^{-2}$ in this case. The two regimes are sketched in \cref{fig:regimes}.

\paragraph{Ghost penalty assumption}
The discrete inequalities facilitated by the Ghost penalty rely on the following assumption which expresses that every active element can be reached from an element with a sufficiently large intersection with the physical domain by a chain of a bounded number of stabilised facets.
\begin{assumption} \label[assumption]{ass:gp}
 For each $n=1,\dots,N$ and each active element $Q_T^n \in \mathcal{Q}_T^n$ there is an element $Q_{T'}^n \in \mathcal{Q}_T^n$ with $|Q_{T'}^n \cap Q^n| \simeq |Q_{T'}^n|$ such that $T$ and $T'$ are connected by a path of at most $M$ facets in $\mathcal{F}^n_{\text{gp}}$, with $M$ independent of $h$, $\Delta t$ and $n$. Moreover, each element $T'$ is the end of at most $M'$ such paths with $M'$ independent of $h$, $\Delta t$ and $n$.
\end{assumption}
This assumption is standard in the analysis of Ghost penalty methods, cf. \cite{burman15, Burman_Hansbo_Larson_Zahedi_2025}; in the space-time setting it has been discussed in detail in \cite[Ass. A.1 and Rem. 4]{preuss18} and \cite[Ass. 4.4]{heimann2025discretizationerroranalysishigh} (for the case of discrete geometries), see also \cite{FH_PHD_2025}. With the choice \cref{eq:def_gp_facets}, the constant $M$ is essentially the number of element layers that are swept by the boundary $\partial \Omega(t)$ within one time slab, i.e. $M \simeq 1 + v_\partial \Delta t / h$, where $v_\partial$ denotes the maximal normal velocity of the boundary. Hence, \cref{ass:gp} is a CFL-type condition on the motion of the \emph{boundary}, $v_\partial \Delta t \lesssim h$, and not on the velocity in the interior of the domain. In Case C it is implied by \cref{eq:ass_caseC}. In Case D it is an additional restriction only for fast moving boundaries; for a stationary or slowly moving boundary the space-time anisotropy $\Delta t \gg h$ allowed by \cref{eq:ass_caseD} is not affected. Note that the restriction of $\mathcal{F}^n_{\text{gp}}$ to the cut region in \cref{eq:def_gp_facets} is no restriction for the assumption: Elements with $Q_T^n \subseteq \overline{Q^n}$ satisfy it trivially (with $T' = T$), and for the remaining elements the paths can be chosen within the cut region, where all facets are stabilised.
The analysis below only uses \cref{ass:gp} and the fact that all cut active facets are stabilised, so that all results remain valid for any larger set of stabilised facets $\mathcal{F}^n_{\text{gp}}$ for which \cref{ass:gp} holds, e.g. all facets between active elements as in \cite{burman2025transport}. For stabilisation strategies for boundaries that sweep many element layers within one time slab we refer to \cite[Ass. A.1 and Rem. 4]{preuss18}, \cite[Ass. 4.4]{heimann2025discretizationerroranalysishigh}.

\paragraph{Discrete convection field}
In Case C the test function will involve the discrete material derivative. As $\w \cdo \nabla u \notin W_h$ in general, we replace $\w$ by an approximation $\wapp$ which is -- on each active prism $Q_T^n$ -- affine in space and constant in time, e.g. the nodal interpolation of $\w(\cdot, t_{n-1})$. Then $\wapp \cdo \nabla u \in W_h$ for $u \in W_h$ and
\begin{equation}
 \Vert \w - \wapp \Vert_{\infty, \mathcal{Q}_T} \lesssim h + \Delta t \lesssim h,
 \quad
 \Vert \wapp \Vert_{\infty} + \Vert \nabla \wapp \Vert_{\infty} \lesssim \Vert \w \Vert_{W^{1,\infty}}, \quad
 \Vert \wapp|_{T_1} - \mathcal{E}^p \wapp|_{T_2} \Vert_{\infty, \omega_F} \lesssim h, \label{eq:wapp_props}
\end{equation}
for all facet patches $\omega_F = T_1 \cup T_2$, where the first estimate uses $\Delta t \lesssim h$ from \cref{eq:ass_caseC} (recall that $\wapp$ is only used in Case C), the second estimate is the $W^{1,\infty}$-stability of the nodal interpolation and the last estimate follows from the fact that both functions are first order approximations of $\w$ on $\omega_F$. We define the discrete directional derivative
\begin{equation}
 \Lapp u := \begin{cases} \partial_t u & \text{ in Case D}, \\ \partial_t u + \wapp \cdo \nabla u & \text{ in Case C}, \end{cases} \qquad \Lapp u \in W_h \text{ for } u \in W_h. \label{eq:def_Lapp}
\end{equation}

\paragraph{Temporal projection}
Finally, we denote by $\Pi_t$ the $L^2(I_n)$-orthogonal projection onto $\mathcal{P}^{k_t}(I_n)$, applied pointwise in space, i.e. for $g \in L^2(\tilde \Omega \times I_n)$, $\Pi_t g(x,\cdot) \in \mathcal{P}^{k_t}(I_n)$ and $\int_{I_n} (g - \Pi_t g)(x,t) q(t) \, \mathrm{d}t = 0$ for all $q \in \mathcal{P}^{k_t}(I_n)$. $\Pi_t$ is $L^2$-stable, commutes with spatial operations such as $\nabla$, the fixed-orientation facet jump $\delta_F$ and the Ghost penalty jump $\jump{\cdot}_{\omega_F}$ (but not with the upwind jump $\fjump{\cdot}$, whose orientation may change in time; as $|\fjump{\cdot}| = |\delta_F \cdot|$ pointwise, all norm estimates for $\fjump{\cdot}$ below are obtained via $\delta_F$), and for $u \in W_h$, $\varphi \in C^\infty([0,T])$ we have $\Pi_t(\varphi u) \in W_h$. We will need the following superapproximation result for products of a discrete function with a smooth function $\varphi = \varphi(t)$, which is known as the \emph{discrete commutator property}, cf. \cite{BERTOLUZZA19991097}. Exploiting the tensor-product structure of the active prisms, it can be applied for a fixed point in space on the respective time element, cf. \cite{burman2025transport}.
\begin{lemma} \label{lemma_superapprox}
 Let $\varphi \in C^{\infty}([0,T])$ and $v \in W_h$. Then for $n=1,\dots,N$, and $k_t \geq 1$ in \cref{eq_superapprox_dt},
 \begin{align}
  \Vert (\mathrm{Id} - \Pi_t) (\varphi v) \Vert_{\mathcal{Q}_T^n} &\leq C(\varphi) \cdot \Delta t \cdot \Vert v \Vert_{\mathcal{Q}_T^n}, \label{eq_superapprox} \\
  \Vert \partial_t (\mathrm{Id} - \Pi_t) (\varphi v) \Vert_{\mathcal{Q}_T^n} &\leq C(\varphi) \cdot \Delta t \cdot \left( \Vert \partial_t v \Vert_{\mathcal{Q}_T^n} + \Vert v \Vert_{\mathcal{Q}_T^n} \right) \leq C(\varphi) \Vert v \Vert_{\mathcal{Q}_T^n}, \label{eq_superapprox_dt}
 \end{align}
 where the last inequality is the inverse inequality \cref{eq:std_inverse} and $\Delta t \leq 1$.
\end{lemma}
For the proof we refer to \cite{BERTOLUZZA19991097, burman2025transport}. For $k_t = 0$ the right hand side bound $C(\varphi) \Vert v \Vert_{\mathcal{Q}_T^n}$, which is all we use below, follows directly from $\partial_t (\mathrm{Id} - \Pi_t)(\varphi v) = \varphi' v$. As the estimates are obtained pointwise in space, they hold with the norms over $\mathcal{Q}_T^n$ replaced on both sides by norms over $S \times I_n$ for any set $S$ on which $v$ is a polynomial in time. In particular, this applies on the full facets $S = F$ to the facet jump with a fixed orientation, $\delta_F v$, and on the facet patches $S = \omega_F$ to the Ghost penalty jump $\jump{v}_{\omega_F}$; both are polynomials in time and commute with $\Pi_t$ and $\varphi$. As $|\fjump{\cdot}| = |\delta_F \cdot|$ pointwise, we obtain for the facet jump and the Ghost penalty jump
\begin{subequations} \label{eq_superapprox_facets}
\begin{align}
 \Vert \fjump{(\mathrm{Id} - \Pi_t) (\varphi v)} \Vert_{F \times I_n} &\leq C(\varphi) \cdot \Delta t \cdot \Vert \fjump{v} \Vert_{F \times I_n}, \label{eq_superapprox_F} \\
 \Vert \jump{(\mathrm{Id} - \Pi_t) (\varphi v)}_{\omega_F} \Vert_{\omega_F \times I_n} &\leq C(\varphi) \cdot \Delta t \cdot \Vert \jump{v}_{\omega_F} \Vert_{\omega_F \times I_n}. \label{eq_superapprox_patch}
\end{align}
\end{subequations}

\subsection{Discrete norms} \label{ssec:norms}
We introduce the following discretisation related norms
\begin{subequations} \label{eq:norms}
\begin{align}
 \tripplejump{u}^2 &:= \sum_{n=1}^{N} \left( \Vert \fjump{u}\Vert_{\Omega^{n-1}}^2
  + \| |\w \cdo \n |^{\frac12} \fjump{u} \|_{\mathcal{Q}_F^n \cap Q}^2  \right)  + \Vert u_- \Vert_{\Omega^N}^2 \label{eq:upwnorm}\\
  & = \Vert |\w^\ast \cdo \n^\ast |^{\frac12} \fjump{u} \Vert_{\mathcal{Q}_F^{\text{st}} \cap Q}^2  + \Vert u_- \Vert_{\Omega^N}^2, \tag{\ref{eq:upwnorm}'}
  \\
  \|u\|_c^2 &:= \Vert u \Vert_Q^2 + \alpha \vert u\vert_{1,h}^2
  + \tripplejump{u}^2, ~~\text{ with }~~  \vert  u\vert_{1,h}^2 :=  \| \nabla u \|_{Q}^2
     +  \| h^{-\frac12}\fjump{u} \|_{\mathcal{Q}_F \cap Q}^2,
   \label{eq_def_coerc_norm}
  \\
  \|u\|_{c,j}^2 &:= \|u\|_c^2 + |u|_J^2, ~~\text{ with }~~ |u|_J^2 := J(u,u), \\
  \tripplenorm{u}^2 &:= \Delta t \| \mathcal{L} u \|_{Q}^2 + \Vert u \Vert_c^2, \qquad
  \tripplenorm{u}^2_j := \tripplenorm{u}^2 + |u|_J^2 = \Delta t \| \mathcal{L} u \|_{Q}^2 + \Vert u \Vert_{c,j}^2. \label{eq_def_tripplenorm}
\end{align}
\end{subequations}
The norm $\Vert \cdot \Vert_c$ will be controlled by symmetric testing, cf. \cref{ssec:coercivity}. It consists of the space-time upwind (semi-)norm $\tripplejump{\cdot}$, the interior penalty norm of the diffusion and the $L^2(Q)$ norm, which is obtained by testing with a temporally weighted version of $u$; $|\cdot|_J$ adds the control of the Ghost penalty. As $\Vert \cdot \Vert_c$ no longer controls any derivative of $u$ for $\alpha \to 0$, the norm $\tripplenorm{\cdot}$ adds the directional derivative $\mathcal{L} u$, i.e. $\partial_t u$ in Case D and the material derivative $\w^\ast \cdo \nabla^\ast u$ in Case C. This control is obtained in \cref{ssec:infsup}.

For the continuity estimate, cf. \cref{ssec:continuity}, we will need the following stronger norm
\begin{equation}
 \tripplenorm{u}_\ast^2 := \tripplenorm{u}^2 + \Delta t^{-1} \Vert u \Vert_Q^2 + \sum_{n=1}^N \Vert u_- \Vert_{\Omega^n}^2 + \frac{h}{\Delta t} \Vert u_- \Vert_{\mathcal{Q}_F \cap Q}^2 + \alpha \Vert h^{\frac12} \avg{\nabla u \cdo \n} \Vert_{\mathcal{Q}_F \cap Q}^2, \label{eq_def_tripplenorm_ast}
\end{equation}
which is well-defined for $u \in W_h + H^2(Q)$. Note that only the upwind traces $u_-$ enter; in the interpolation estimates of \cref{sect_a_pr_est} we will simply bound them by the direction-independent quantities involving both one-sided traces $u_\pm$.

\subsection{Ghost penalty related estimates} \label{ssec:gp}
The numerical analysis relies on discrete inequalities, which in the cut case are facilitated by the Ghost penalty stabilisation \cite{burman15, Burman_Hansbo_Larson_Zahedi_2025}. As it relates to the time-dependent direct Ghost penalty, in the works \cite{preuss18, heimann20, heimann2025discretizationerroranalysishigh, burman2025transport} such discrete inequalities were established. We recall them here in our notation. All of them rely on \cref{ass:gp}.

\begin{lemma} \label{lemma_gp_E_IpGP}
 It holds for every $u_h \in W_h$, $n=1,\dots,N$,
\begin{align}
 \|u_h\|^2_{\mathcal{Q}_T^n} &\lesssim  j^n_h(u_h,u_h) + \| u_h \|_{Q^n}^2 \quad \textnormal{and hence} \quad \|u_h\|^2_{\mathcal{Q}_T} \lesssim c_J^{-1} |u_h|_J^2 + \| u_h \|_{Q}^2. \label{eq1_gp_E_IpGP}
\end{align}
\end{lemma}
\begin{proof}
 We refer the reader to \cite{preuss18, heimann2025discretizationerroranalysishigh} and note that in contrast to these works no scaling in $h$ is contained in $j^n_h$ here. The key observation is that for a facet patch $\omega_F = T_1 \cup T_2$ and a polynomial $u_h$ on $\omega_F$ we have $\Vert u_h \Vert_{T_1 \times I_n}^2 \lesssim \Vert \jump{u_h}_{\omega_F} \Vert_{T_1 \times I_n}^2 + \Vert u_h \Vert_{T_2 \times I_n}^2$ and $\Vert u_h \Vert_{T_2 \times I_n}^2 \lesssim \Vert u_h \Vert_{Q_{T_2}^n \cap Q^n}^2$ if $|Q_{T_2}^n \cap Q^n| \simeq |Q_{T_2}^n|$, by norm equivalences for polynomials (the latter being a norm equivalence between a subset of comparable measure and the full prism, cf. \cite{preuss18}). Following the paths of \cref{ass:gp} yields the claim.
\end{proof}
\begin{lemma} \label{lemma_bound_tn_ag_temporal_vol}
For any $u \in W_h$, $n=1,\dots,N$, it holds
\begin{align}
 \Vert u_+ \Vert_{\Omega^{n-1}}^2 + \Vert u_- \Vert_{\Omega^{n}}^2 &\lesssim \frac{1}{\Delta t} \|u\|_{\mathcal{Q}_T^{n}}^2
 \lesssim \frac{1}{\Delta t} \left( j^n_h(u,u) + \|u\|_{Q^{n}}^2 \right). \label{eq_bound_tn_ag_temporal_vol}
\end{align}
\end{lemma}
\begin{lemma} \label{lemma_temporal_and_spatial_inverse_ineq_gp}
 For any $u_h \in W_h$, $n=1,\dots,N$, it holds
\begin{align}
 \Vert \Delta t \partial_t u_h \Vert_{Q^n}^2 + \Vert h \nabla u_h \Vert_{Q^n}^2 &\lesssim \|u_h\|^2_{\mathcal{Q}_T^n} \lesssim j^n_h(u_h,u_h) + \| u_h \|_{Q^n}^2. \label{eq_temporal_inverse_ineq_gp}
\end{align}
\end{lemma}
\begin{proof}
 Both lemmas follow from the respective discrete trace and inverse inequalities on the full prisms, \cref{eq:std_trace} and \cref{eq:std_inverse}, summation over all active elements and an application of \cref{lemma_gp_E_IpGP}. For \cref{lemma_bound_tn_ag_temporal_vol} we note that $\Omega^{n-1} \cup \Omega^n$ is covered by the active elements of $\mathcal{Q}_T^n$.
\end{proof}

\begin{lemma}\label{gp_inv_ineq}
 For any $u\in W_h$, it holds
\begin{align}
 j_h^n( h \nabla u, h \nabla u) + j_h^n( \Delta t \, \partial_t  u, \Delta t \, \partial_t  u) &\lesssim j_h^n(u,u), \\ \text{and hence} \quad | h \nabla u|_J^2 + | \Delta t \, \partial_t  u |_J^2 \lesssim | u |_J^2. \label{eq_gp_inv_ineq_dt_grad}
\end{align}
\end{lemma}
\begin{proof}
 As $\mathcal{E}^p$ commutes with $\nabla$ and $\partial_t$, we have $\jump{\nabla u}_{\omega_F} = \nabla \jump{u}_{\omega_F}$ and $\jump{\partial_t u}_{\omega_F} = \partial_t \jump{u}_{\omega_F}$ on each $T_i \times I_n$, where $\jump{u}_{\omega_F}|_{T_i}$ is a polynomial. The inverse inequalities \cref{eq:std_inverse} yield the claim.
\end{proof}

A new ingredient of the DG-in-space setting is that the Ghost penalty also controls the DG jumps across the \emph{full} stabilised facets $F \times I_n$, $F \in \mathcal{F}^n_{\text{gp}}$, which is not the case for the jump terms $\Vert \fjump{u} \Vert_{\mathcal{Q}_F \cap Q}$ that appear in the bilinear form. On the interior facets $F \in \mathcal{F}^n_{\text{int}}$ no such control is available, but none is needed either, as $F \times I_n \subseteq \overline{Q^n}$ there, cf. \cref{eq:def_int_facets}.
\begin{lemma} \label{lemma_gp_facet_jumps}
 For any $u \in W_h$ it holds
 \begin{equation}
  \alpha \Vert h^{-\frac12} \fjump{u} \Vert_{\mathcal{Q}_F^{\text{gp}}}^2 + \Vert \w \Vert_\infty \Vert \fjump{u} \Vert_{\mathcal{Q}_F^{\text{gp}}}^2 \lesssim \left( \frac{\alpha}{h^2} + \frac{\Vert \w \Vert_\infty}{h} \right) \sum_{n=1}^N j_h^n(u,u) = \gamma_J^{-1} |u|_J^2. \label{eq_gp_facet_jumps}
 \end{equation}
\end{lemma}
\begin{proof}
 Let $F = \overline{T_1} \cap \overline{T_2} \in \mathcal{F}_h^n \cap \mathcal{F}^n_{\text{gp}}$ be a stabilised active facet. As $\mathcal{E}^p(u|_{T_2}) = u|_{T_2}$ on $F$, we have $\fjump{u}|_{F \times I_n} = \pm \jump{u}_{\omega_F}|_{T_1} |_{F \times I_n}$, where $\jump{u}_{\omega_F}|_{T_1}$ is a polynomial on $T_1 \times I_n$. The discrete trace inequality \cref{eq:std_trace} yields $\Vert \fjump{u} \Vert_{F \times I_n}^2 \lesssim h^{-1} \Vert \jump{u}_{\omega_F} \Vert_{T_1 \times I_n}^2 \leq h^{-1} j_F^n(u,u)$. Summation over all stabilised facets and time slabs yields the claim.
\end{proof}

Next, we collect the properties of the discrete directional derivative $\Lapp$, cf. \cref{eq:def_Lapp}.
\begin{lemma} \label{lemma_Lapp}
 For any $u \in W_h$ it holds (in both cases)
 \begin{subequations}
 \begin{align}
  \mathcal{L} u = \Lapp u \ \text{(Case D)}, &\qquad \Vert \mathcal{L} u - \Lapp u \Vert_{Q} \lesssim h \Vert \nabla u \Vert_Q \ \text{(Case C)}, \label{eq_Lapp_diff} \\
  \Vert \Lapp u \Vert_{\mathcal{Q}_T^n} \lesssim \Delta t^{-1} \Vert u \Vert_{\mathcal{Q}_T^n}, &\qquad
  \Vert \nabla \Lapp u \Vert_{\mathcal{Q}_T^n} \lesssim \Delta t^{-1} \Vert \nabla u \Vert_{\mathcal{Q}_T^n}, \label{eq_Lapp_inverse} \\
  j_h^n(\Lapp u, \Lapp u) &\lesssim \Delta t^{-2} j_h^n(u,u) + h^2 \Vert \nabla u \Vert_{\mathcal{Q}_T^n}^2, \label{eq_Lapp_gp}
 \end{align}
 \end{subequations}
 where in Case D the last term in \cref{eq_Lapp_gp} can be dropped.
\end{lemma}
\begin{proof}
 In Case D, $\mathcal{L} u = \Lapp u = \partial_t u$ by definition and \cref{eq_Lapp_inverse}, \cref{eq_Lapp_gp} are the inverse inequalities \cref{eq:std_inverse} and \cref{gp_inv_ineq}.
 In Case C, \cref{eq_Lapp_diff} follows from $\mathcal{L} u - \Lapp u = (\w - \wapp) \cdo \nabla u$ and \cref{eq:wapp_props}. For \cref{eq_Lapp_inverse} we use \cref{eq:std_inverse}, \cref{eq:wapp_props} and the CFL-type condition $\Vert \w \Vert_\infty \Delta t \lesssim h$ from \cref{eq:ass_caseC}:
 $\Vert \wapp \cdo \nabla u \Vert_{Q_T^n} \lesssim \Vert \w \Vert_\infty h^{-1} \Vert u \Vert_{Q_T^n} \lesssim \Delta t^{-1} \Vert u \Vert_{Q_T^n}$ and $\Vert \nabla (\wapp \cdo \nabla u) \Vert_{Q_T^n} \lesssim (\Vert \nabla \wapp \Vert_\infty + \Vert \w \Vert_\infty h^{-1}) \Vert \nabla u \Vert_{Q_T^n} \lesssim \Delta t^{-1} \Vert \nabla u \Vert_{Q_T^n}$.
 For \cref{eq_Lapp_gp} we note that the polynomial extension of a product is the product of the extensions, so that on $T_i \times I_n$, $i=1,2$,
 \begin{equation*}
  \jump{\wapp \cdo \nabla u}_{\omega_F} = \wapp|_{T_i} \cdo \nabla \jump{u}_{\omega_F} + (\wapp|_{T_i} - \mathcal{E}^p \wapp|_{T_{3-i}}) \cdo \mathcal{E}^p (\nabla u |_{T_{3-i}}).
 \end{equation*}
 The first term is bounded by $\Vert \w \Vert_\infty h^{-1} \Vert \jump{u}_{\omega_F} \Vert_{T_i \times I_n} \lesssim \Delta t^{-1} \Vert \jump{u}_{\omega_F} \Vert_{T_i \times I_n}$ by \cref{eq:std_inverse}, the second one by $h \Vert \nabla u \Vert_{T_{3-i} \times I_n}$ by \cref{eq:wapp_props} and the stability of the polynomial extension on shape regular patches. Together with \cref{gp_inv_ineq} for the $\partial_t$ part and summation over the facets (each element belongs to a bounded number of patches) this yields \cref{eq_Lapp_gp}.
\end{proof}

The following consequence is the central estimate for the derivative test function $\Delta t \Lapp u$ in the inf-sup analysis.
\begin{lemma} \label{lemma_Lapp_tripplenorm}
 For any $u \in W_h$ it holds (in both cases)
 \begin{align}
  \Delta t \Vert \Lapp u \Vert_{\mathcal{Q}_T}^2 + \vert \Delta t \Lapp u \vert_J^2 &\lesssim \Delta t \Vert \mathcal{L} u \Vert_Q^2 + \Vert u \Vert_{c,j}^2 = \tripplenorm{u}_j^2, \label{eq_Lapp_tripplenorm} \\
  \alpha \Vert \Delta t \nabla \Lapp u \Vert_{\mathcal{Q}_T}^2 \lesssim \alpha \Vert \nabla u \Vert_{\mathcal{Q}_T}^2 &\lesssim \alpha \Vert \nabla u \Vert_Q^2 + \gamma_J^{-1} |u|_J^2 \lesssim \Vert u \Vert_{c,j}^2. \label{eq_alpha_grad_full}
 \end{align}
\end{lemma}
\begin{proof}
 We start with the first term. By \cref{lemma_gp_E_IpGP}, \cref{eq_Lapp_diff} and \cref{eq_Lapp_gp}
 \begin{align*}
  \Delta t \Vert \Lapp u \Vert_{\mathcal{Q}_T}^2 &\lesssim \Delta t \Vert \Lapp u \Vert_Q^2 + \Delta t \sum_{n} j_h^n(\Lapp u, \Lapp u) \\
  & \lesssim \Delta t \Vert \mathcal{L} u \Vert_Q^2 + \Delta t \Vert \mathcal{L} u - \Lapp u \Vert_Q^2 + \Delta t^{-1} \sum_n j_h^n(u,u) + \Delta t h^2 \Vert \nabla u \Vert_{\mathcal{Q}_T}^2.
 \end{align*}
 By \cref{eq:ass_gp_scaling}, $\Delta t^{-1} \sum_n j_h^n(u,u) \lesssim c_J^{-1} \Delta t^{-1} |u|_J^2 \lesssim |u|_J^2$.
 Further, by the inverse inequality \cref{eq:std_inverse} on the full prisms, \cref{lemma_gp_E_IpGP}, \cref{eq:ass_gp_scaling} and $\Delta t \leq 1$, $\Delta t h^2 \Vert \nabla u \Vert_{\mathcal{Q}_T}^2 \lesssim \Delta t \Vert u \Vert_{\mathcal{Q}_T}^2 \lesssim \Delta t \Vert u \Vert_Q^2 + \Delta t c_J^{-1} |u|_J^2 \lesssim \Vert u \Vert_{c,j}^2$.
 The difference term vanishes in Case D, and in Case C we have $\Delta t \Vert \mathcal{L} u - \Lapp u \Vert_Q^2 \lesssim \Delta t h^2 \Vert \nabla u \Vert_Q^2 \lesssim \Vert u \Vert_{c,j}^2$ as before.

 For the Ghost penalty term we have by \cref{eq_Lapp_gp}
 \begin{align*}
  \vert \Delta t \Lapp u \vert_J^2 = c_J \Delta t^2 \sum_n j_h^n(\Lapp u, \Lapp u) \lesssim |u|_J^2 + c_J \Delta t^2 h^2 \Vert \nabla u \Vert_{\mathcal{Q}_T}^2,
 \end{align*}
 where the second term only arises in Case C, cf. \cref{eq_Lapp_gp}. By the definition of $c_J$, \cref{lemma_gp_E_IpGP} applied to $\nabla u \in [W_h]^d$, \cref{gp_inv_ineq} and \cref{eq:std_inverse},
 \begin{align*}
  c_J \Delta t^2 h^2 \Vert \nabla u \Vert_{\mathcal{Q}_T}^2 &\lesssim \gamma_J \Delta t^2 \alpha \Vert \nabla u \Vert_{\mathcal{Q}_T}^2 + \gamma_J \Vert \w \Vert_\infty h \Delta t^2 \Vert \nabla u \Vert_{\mathcal{Q}_T}^2 \\
  & \lesssim \Delta t^2 \Big( \alpha \Vert \nabla u \Vert_{Q}^2 + \alpha \sum_n j_h^n(\nabla u, \nabla u) \Big) + \Delta t \Vert u \Vert_{\mathcal{Q}_T}^2 \\
  & \lesssim \alpha \Vert \nabla u \Vert_Q^2 + |u|_J^2 + \Vert u \Vert_{c,j}^2,
 \end{align*}
 where for the second term we used $\Vert \w \Vert_\infty h \Delta t^2 \Vert \nabla u \Vert_{\mathcal{Q}_T}^2 \lesssim \Vert \w \Vert_\infty \Delta t^2 h^{-1} \Vert u \Vert_{\mathcal{Q}_T}^2 \lesssim \Delta t \Vert u \Vert_{\mathcal{Q}_T}^2$ by \cref{eq:std_inverse} and the CFL-type condition of Case C, and $\Delta t \leq 1$ for the first term.
 Finally, \cref{eq_alpha_grad_full} follows from \cref{eq_Lapp_inverse}, \cref{lemma_gp_E_IpGP} applied to $\nabla u \in [W_h]^d$ and \cref{gp_inv_ineq}:
 \begin{equation*}
  \alpha \Vert \Delta t \nabla \Lapp u \Vert_{\mathcal{Q}_T}^2 \lesssim \alpha \Vert \nabla u \Vert_{\mathcal{Q}_T}^2 \lesssim \alpha \Vert \nabla u \Vert_Q^2 + \alpha h^{-2} \sum_n j_h^n(u,u) \lesssim \alpha \Vert \nabla u \Vert_Q^2 + \gamma_J^{-1} |u|_J^2.
 \end{equation*}
\end{proof}
Note that \cref{eq_alpha_grad_full}, in particular $\alpha \Vert \nabla u \Vert_{\mathcal{Q}_T}^2 \lesssim \alpha \Vert \nabla u \Vert_Q^2 + \gamma_J^{-1} |u|_J^2$, will be used repeatedly below. It is the reason for the scaling $\alpha h^{-2}$ in $c_J$.

We close this section with two estimates for DG jumps that do not rely on the Ghost penalty on the interior facets. The first one shows that jumps weighted with $\Vert \w \Vert_\infty$ (rather than with $|\w \cdo \n|$ as in $\tripplejump{\cdot}$) are always controlled by other quantities.
\begin{lemma} \label{lemma_w_jumps}
 For any $v$ that is a piecewise polynomial in space on the active mesh, i.e. $v|_{Q_T^n} \in \mathcal{P}^{k_s}(T) \otimes L^2(I_n)$ for all $Q_T^n \in \mathcal{Q}_T$ (in particular for $v \in W_h$ and for $v = (\Pi_t - \mathrm{Id})(\varphi u)$ with $u \in W_h$), it holds
 \begin{align}
  \Vert \w \Vert_\infty \Vert \fjump{v} \Vert_{\mathcal{Q}_F \cap Q}^2 &\lesssim \alpha \Vert h^{-\frac12} \fjump{v} \Vert_{\mathcal{Q}_F \cap Q}^2 \quad \text{(Case D)}, \label{eq_w_jumps} \\
  \Vert \w \Vert_\infty \Vert \fjump{v} \Vert_{\mathcal{Q}_F}^2 &\lesssim \Delta t^{-1} \Vert v \Vert_{\mathcal{Q}_T}^2 \quad \text{(Case C)}. \nonumber
 \end{align}
\end{lemma}
\begin{proof}
 In Case D we have $\Vert \w \Vert_\infty \lesssim \alpha h^{-1}$ by \cref{eq:ass_caseD}. In Case C, the discrete trace inequality \cref{eq:std_trace} in space (applied pointwise in time) gives $\Vert \fjump{v} \Vert_{F \times I_n}^2 \lesssim h^{-1} \Vert v \Vert_{\omega_F \times I_n}^2$ for every active facet, and $\Vert \w \Vert_\infty h^{-1} \lesssim \Delta t^{-1}$ by the CFL-type condition \cref{eq:ass_caseC}.
\end{proof}
The second one is the jump estimate for the derivative test function $\Delta t \Lapp u$. On the stabilised facets it follows from the Ghost penalty, on the interior facets we exploit that $F \times I_n \subseteq \overline{Q^n}$ is a tensor-product domain.
\begin{lemma} \label{lemma_y_jumps}
 For any $u \in W_h$ and $y := \Delta t \Lapp u$ it holds (in both cases)
 \begin{equation}
  \alpha \Vert h^{-\frac12} \fjump{y} \Vert_{\mathcal{Q}_F \cap Q}^2 + \Vert \w \Vert_\infty \Vert \fjump{y} \Vert_{\mathcal{Q}_F \cap Q}^2 \lesssim \tripplenorm{u}_j^2. \label{eq_y_jumps}
 \end{equation}
\end{lemma}
\begin{proof}
 By \cref{lemma_w_jumps} the second term is bounded by the first one in Case D and by $\Delta t \Vert \Lapp u \Vert_{\mathcal{Q}_T}^2 \lesssim \tripplenorm{u}_j^2$ in Case C, cf. \cref{eq_Lapp_tripplenorm}. For the first term we split $\mathcal{Q}_F \cap Q \subseteq \mathcal{Q}_F^{\text{gp}} \cup \mathcal{Q}_F^{\text{int}}$. On the stabilised facets, \cref{lemma_gp_facet_jumps} and \cref{eq_Lapp_tripplenorm} yield $\alpha \Vert h^{-\frac12} \fjump{y} \Vert_{\mathcal{Q}_F^{\text{gp}}}^2 \lesssim \gamma_J^{-1} |y|_J^2 \lesssim \tripplenorm{u}_j^2$.
 Now let $F \in \mathcal{F}^n_{\text{int}}$ be an interior facet. Then $F \times I_n \subseteq \overline{Q^n}$, so that $\Vert \fjump{u} \Vert_{F \times I_n} = \Vert \fjump{u} \Vert_{F \times I_n \cap Q}$ is controlled by $\Vert u \Vert_c$, and $\delta_F u$ is a polynomial in time on the tensor-product domain $F \times I_n$. In Case D we have $|\fjump{y}| = \Delta t |\partial_t \delta_F u|$ (as $\delta_F$ commutes with $\partial_t$) and the inverse inequality in time, applied for each fixed $x \in F$, yields
 \begin{equation*}
  \Vert \fjump{y} \Vert_{F \times I_n}^2 = \Delta t^2 \Vert \partial_t \delta_F u \Vert_{F \times I_n}^2 \lesssim \Vert \delta_F u \Vert_{F \times I_n}^2 = \Vert \fjump{u} \Vert_{F \times I_n \cap Q}^2.
 \end{equation*}
 In Case C, $|\fjump{y}| \leq \Delta t |\partial_t \delta_F u| + \Delta t |\delta_F (\wapp \cdo \nabla u)|$, where the first part is treated as in Case D. For the second part, the trace inequality \cref{eq:std_trace} applied to both sides of $F$, \cref{eq:wapp_props} and the CFL-type condition \cref{eq:ass_caseC} give
 \begin{equation*}
  \alpha h^{-1} \Delta t^2 \Vert \delta_F (\wapp \cdo \nabla u) \Vert_{F \times I_n}^2 \lesssim \alpha h^{-2} \Delta t^2 \Vert \w \Vert_\infty^2 \Vert \nabla u \Vert_{\omega_F \times I_n}^2 \lesssim \alpha \Vert \nabla u \Vert_{\omega_F \times I_n}^2.
 \end{equation*}
 Summation over the interior facets and time slabs (each element belongs to a bounded number of facet patches) and \cref{eq_alpha_grad_full} yield
 \begin{equation*}
  \alpha \Vert h^{-\frac12} \fjump{y} \Vert_{\mathcal{Q}_F^{\text{int}}}^2 \lesssim \alpha \Vert h^{-\frac12} \fjump{u} \Vert_{\mathcal{Q}_F \cap Q}^2 + \alpha \Vert \nabla u \Vert_{\mathcal{Q}_T}^2 \lesssim \Vert u \Vert_{c,j}^2.
 \end{equation*}
\end{proof}

\subsection{Partial integration and symmetric testing} \label{ssec:coercivity}
In this section, we investigate the contribution from considering $B_h(u,u)$ and $B_h(u,\Pi_t(e^{-\sigma t} u))$, $\sigma > 0$, in the stability argument. We start with the partial integration identity for the convective part, which is the discrete analogue of the fact that the material derivative is skew-symmetric for a divergence free convection field and a domain that moves with the flow.
\begin{lemma}[Mass conserving form] \label{lemma_B_Bmc}
 For all $u, v$ with $u, v \in H^1(Q_T^n \cap Q)$ for every active prism $Q_T^n \in \mathcal{Q}_T$ (in particular for $u, v \in W_h + H^1(Q)$ and for $v = \varphi u$ with $u \in W_h$ and $\varphi \in C^1([0,T])$) it holds
 \begin{align}
  (C + U_h)(u,v) = B_{mc}(u,v) :={}& - (u, \w^\ast \cdo \nabla^\ast v)_Q \nonumber \\
  & - (|\w^\ast \cdo \n^\ast| u_-, \fjump{v})_{\mathcal{Q}_F^{\text{st}} \cap Q} + (u_-, v_-)_{\Omega^N}, \label{eq_B_Bmc}
 \end{align}
 with the convention $u_- := 0$ on $\Omega^0$. In particular, the right hand side does not involve derivatives of $u$, and
 \begin{align}
  (C + U_h)(u,u) = \tfrac12 \tripplejump{u}^2. \label{eq_CU_symm}
 \end{align}
\end{lemma}
\begin{proof}
 We apply partial integration in space-time on each cut space-time prism $Q_T^n \cap Q$, $Q_T^n \in \mathcal{Q}_T^n$, and use $\operatorname{div}^\ast \w^\ast = \operatorname{div} \w = 0$:
 \begin{equation*}
  (\w^\ast \cdo \nabla^\ast u, v)_{Q_T^n \cap Q} = - (u, \w^\ast \cdo \nabla^\ast v)_{Q_T^n \cap Q} + \int_{\partial (Q_T^n \cap Q)} u v \, \w^\ast \cdo \n^\ast_{\partial} \, \mathrm{d}s,
 \end{equation*}
 with $\n^\ast_\partial$ the outer unit normal. The boundary $\partial (Q_T^n \cap Q)$ consists of parts of the lateral boundary of $Q$, where $\w^\ast \cdo \n^\ast_\partial = 0$ as the domain moves with the velocity $\w$, cf. \cref{eq:domain_motion}, the top and bottom parts $(T \cap \Omega^{n}) \times \{t_n\}$ and $(T \cap \Omega^{n-1}) \times \{t_{n-1}\}$, where $\w^\ast \cdo \n^\ast_\partial = \pm 1$, and parts of the spatial facets $(F \times I_n) \cap Q$. Summing over all elements, each interior facet appears twice, with $\n_\partial^\ast = \n^\ast$ from the upwind element (with values $u_-, v_-$) and $\n_\partial^\ast = -\n^\ast$ from the downwind element (with values $u_+, v_+$). Hence,
 \begin{align*}
  C^n(u,v) &= - (u, \w^\ast \cdo \nabla^\ast v)_{Q^n} + (u_-, v_-)_{\Omega^n} - (u_+, v_+)_{\Omega^{n-1}} \\ & \quad + (|\w \cdo \n|, u_- v_- - u_+ v_+)_{\mathcal{Q}_F^n \cap Q}.
 \end{align*}
 Adding $U_h^n(u,v)$ and using $u_- v_- - u_+ v_+ + (u_+ - u_-) v_+ = - u_- \fjump{v}$ pointwise on the spatial facets and on $\Omega^{n-1}$, we obtain
 \begin{align*}
  (C^n + U_h^n)(u,v) &= - (u, \w^\ast \cdo \nabla^\ast v)_{Q^n} - (|\w \cdo \n| u_-, \fjump{v})_{\mathcal{Q}_F^n \cap Q} \\ & \quad + (u_-, v_-)_{\Omega^n} - (u_-, v_+)_{\Omega^{n-1}}.
 \end{align*}
 Summation over $n$ and rearranging $\sum_{n=1}^N (u_-, v_-)_{\Omega^n} - (u_-, v_+)_{\Omega^{n-1}} = (u_-,v_-)_{\Omega^N} - \sum_{n=1}^{N-1} (u_-, \fjump{v})_{\Omega^n}$ (recall $u_- = 0$ on $\Omega^0$) yields \cref{eq_B_Bmc}. For \cref{eq_CU_symm} we add the definition and \cref{eq_B_Bmc} for $u=v$ and use $\fjump{u} u_+ - u_- \fjump{u} = \fjump{u}^2$ pointwise on $\mathcal{Q}_F^{\text{st}} \cap Q$ (on $\Omega^0$ this reads $u_+^2 = \fjump{u}^2$).
\end{proof}

For the diffusive part we have the following standard interior penalty coercivity result, which we state in a slightly more general, temporally weighted form.
\begin{lemma} \label{lemma_IP_coercivity}
 Let $\varphi \in C([0,T])$, $\varphi > 0$, with $\sup_{I_n} \varphi \leq C_\varphi \inf_{I_n} \varphi$ for $n=1,\dots,N$. There is $\lambda_0 > 0$ depending on $C_\varphi$, $\gamma_J$, the shape regularity, the polynomial degrees and the constants in \cref{ass:gp} such that for all $\lambda \geq \lambda_0$ and all $u \in W_h$ there holds
 \begin{equation}
  (A + N_h)(u, \varphi u) \geq \frac{\alpha}{2} \left( \Vert \varphi^{\frac12} \nabla u \Vert_Q^2 + \Vert \varphi^{\frac12} h^{-\frac12} \fjump{u} \Vert_{\mathcal{Q}_F \cap Q}^2 \right) - \frac14 J(u, \varphi u). \label{eq_IP_coercivity}
 \end{equation}
\end{lemma}
\begin{proof}
 As $\varphi$ depends on $t$ only, it commutes with all spatial operations, $\nabla(\varphi u) = \varphi \nabla u$, $\fjump{\varphi u} = \varphi \fjump{u}$, so that
 \begin{equation*}
  (A + N_h)(u, \varphi u) = \alpha \Vert \varphi^{\frac12} \nabla u \Vert_Q^2 + 2 (\avg{\alpha \nabla u \cdo \n}, \varphi \fjump{u})_{\mathcal{Q}_F \cap Q} + \alpha \lambda \Vert \varphi^{\frac12} h^{-\frac12} \fjump{u} \Vert_{\mathcal{Q}_F \cap Q}^2.
 \end{equation*}
 With Cauchy-Schwarz and Young's inequality with $\epsilon > 0$ the mixed term is bounded by
 $\frac{\alpha}{\epsilon} \Vert \varphi^{\frac12} h^{\frac12} \avg{\nabla u \cdo \n} \Vert_{\mathcal{Q}_F}^2 + \epsilon \alpha \Vert \varphi^{\frac12} h^{-\frac12} \fjump{u} \Vert_{\mathcal{Q}_F \cap Q}^2$.
 For the first term, we set $m_n := \inf_{I_n} \varphi$ and $M_n := \sup_{I_n} \varphi \leq C_\varphi m_n$ and apply, slab by slab, the discrete trace inequality \cref{eq:std_trace}, \cref{lemma_gp_E_IpGP} applied to $\nabla u \in [W_h]^d$ and \cref{gp_inv_ineq}:
 \begin{align*}
  \Vert \varphi^{\frac12} h^{\frac12} \avg{\nabla u \cdo \n} \Vert_{\mathcal{Q}_F^n}^2 &\leq C_1 M_n \Vert \nabla u \Vert_{\mathcal{Q}_T^n}^2 \leq C_2 M_n \left( \Vert \nabla u \Vert_{Q^n}^2 + h^{-2} j_h^n(u,u) \right) \\
  &\leq C_2 C_\varphi \left( \Vert \varphi^{\frac12} \nabla u \Vert_{Q^n}^2 + h^{-2} m_n j_h^n(u,u) \right).
 \end{align*}
 As $\varphi$ enters the Ghost penalty as a positive weight in the time integral, $J(u, \varphi u) \geq c_J \sum_n m_n j_h^n(u,u)$, so that $\alpha h^{-2} \sum_n m_n j_h^n(u,u) \leq \gamma_J^{-1} J(u, \varphi u)$ by the definition of $c_J$. Choosing $\epsilon = \max\{ \frac12, 2 C_2 C_\varphi, 4 C_2 C_\varphi \gamma_J^{-1} \}$ and $\lambda_0 = 2 \epsilon$ yields the claim: the gradient term keeps the factor $1 - C_2 C_\varphi / \epsilon \geq \frac12$, the jump term the factor $\lambda - \epsilon \geq \epsilon \geq \frac12$, and the Ghost penalty term has the factor $C_2 C_\varphi \gamma_J^{-1} / \epsilon \leq \frac14$.
\end{proof}
For the weights $\varphi \equiv 1$ and $\varphi = e^{-\sigma t}$ with $\sigma \Delta t \leq 1$ used below, we have $C_\varphi \leq e$. From now on, we assume $\lambda \geq \lambda_0$, where $\lambda_0$ is chosen for $C_\varphi = e$. Combining the two previous lemmas yields the following result for the symmetric testing.
\begin{proposition} \label{prop_symm_testing}
 For all $u_h \in W_h$, it holds
 \begin{equation}
  (B_h + J)(u_h,u_h) \geq \frac12 \tripplejump{u_h}^2 + \frac{\alpha}{2} | u_h|_{1,h}^2 + \frac34 |u_h|_J^2. \label{eq_symm_testing}
 \end{equation}
\end{proposition}
\begin{proof}
 Apply \cref{eq_CU_symm} and \cref{eq_IP_coercivity} with $\varphi \equiv 1$.
\end{proof}
Unlike in the continuous-in-space case \cite{heimann2025discretizationerroranalysishigh}, the Ghost penalty is thus needed already for coercivity, namely to bound the normal derivatives in the interior penalty consistency terms on cut elements.

The estimate \cref{eq_symm_testing} does not yield control on the $L^2(Q)$-norm of $u_h$ (as $\operatorname{div} \w = 0$). To obtain this control we follow \cite{burman2025transport} and test with a discrete version of $e^{-\sigma t} u_h$ for some fixed $\sigma > 0$. In order to keep the dependence on the final time transparent we will keep $\sigma$ as a parameter in the following; the natural choice is $\sigma \simeq 1$ (or $\sigma \simeq 1/T$).
\begin{lemma} \label{lemma_weighted_symm_testing}
 Let $\sigma > 0$ with $\sigma \Delta t \leq 1$ and $u_w := \Pi_t(e^{-\sigma t} u) \in W_h$ for $u \in W_h$. There are constants $c_\sigma, C > 0$, $c_\sigma \simeq e^{-\sigma T} \min\{1, \sigma\}$, independent of $h$, $\Delta t$, $\alpha$, such that for all $u \in W_h$
 \begin{equation}
  (B_h + J)(u, u_w) \geq c_\sigma \Vert u \Vert_{c,j}^2 - C \Delta t^{\frac12} \tripplenorm{u}_j^2. \label{eq_weighted_symm_testing}
 \end{equation}
\end{lemma}
\begin{proof}
 We split $u_w = \varphi u + r$ with $\varphi := e^{-\sigma t}$ and $r := (\Pi_t - \mathrm{Id})(\varphi u)$.

 \underline{Main part:}
 By \cref{eq_B_Bmc} applied to $v = \varphi u$, noting $\w^\ast \cdo \nabla^\ast (\varphi u) = \varphi \w^\ast \cdo \nabla^\ast u + \varphi' u$ and $\fjump{\varphi u} = \varphi \fjump{u}$, and adding the definition of $(C+U_h)(u, \varphi u)$ we obtain as in the proof of \cref{eq_CU_symm}
 \begin{align*}
  (C + U_h)(u, \varphi u) &= \frac12 \tripplejump{\varphi^{\frac12} u}^2 - \frac12 (\varphi' u, u)_Q = \frac12 \tripplejump{\varphi^{\frac12} u}^2 + \frac{\sigma}{2} \Vert \varphi^{\frac12} u \Vert_Q^2 \\
  &\geq \frac{e^{-\sigma T}}{2} \left( \tripplejump{u}^2 + \sigma \Vert u \Vert_Q^2 \right).
 \end{align*}
 Moreover, by \cref{lemma_IP_coercivity} and $J(u, \varphi u) \geq e^{-\sigma T} J(u,u)$ (as $\mathcal{E}^p$ acts pointwise in time, $\jump{\varphi u}_{\omega_F} = \varphi \jump{u}_{\omega_F}$, so that the weight enters the Ghost penalty as a positive weight in the time integral),
 \begin{align*}
  (A + N_h + J)(u, \varphi u) \geq \frac{\alpha}{2} e^{-\sigma T} |u|_{1,h}^2 + \frac34 J(u, \varphi u) \geq e^{-\sigma T} \left( \frac{\alpha}{2} |u|_{1,h}^2 + \frac34 |u|_J^2 \right).
 \end{align*}

 \underline{Remainder:}
 We collect the properties of $r$ that follow from \cref{lemma_superapprox} and \cref{eq_superapprox_facets} (applied to $u$, $\nabla u$, $\delta_F u$ and $\jump{u}_{\omega_F}$, as $\Pi_t$ commutes with these spatial operations, and $|\fjump{r}| = |\delta_F r| = |(\Pi_t - \mathrm{Id})(\varphi \delta_F u)|$ pointwise), the inverse inequality \cref{eq:std_inverse} and the trace inequalities \cref{eq:std_trace} (in space; $r$ is a polynomial in space) and \cref{eq:cont_trace} (in time):
 \begin{subequations} \label{eq_r_props}
 \begin{align}
  \Vert r \Vert_{\mathcal{Q}_T^n} \lesssim \Delta t \Vert u \Vert_{\mathcal{Q}_T^n}, \quad &\Vert \partial_t r \Vert_{\mathcal{Q}_T^n} \lesssim \Vert u \Vert_{\mathcal{Q}_T^n}, \quad \Vert \nabla r \Vert_{\mathcal{Q}_T^n} \lesssim \Delta t \Vert \nabla u \Vert_{\mathcal{Q}_T^n}, \\
  \Vert \fjump{r} \Vert_{\mathcal{Q}_F^n} \lesssim \Delta t \Vert \fjump{u} \Vert_{\mathcal{Q}_F^n}, \quad &|r|_J \lesssim \Delta t |u|_J, \quad
  \Vert r_\pm \Vert_{\mathcal{Q}_F^n} \lesssim h^{-\frac12} \Delta t \Vert u \Vert_{\mathcal{Q}_T^n}, \label{eq_r_props_b} \\
  \Vert r_+ \Vert_{\Omega^{n-1}}^2 + \Vert r_- \Vert_{\Omega^{n}}^2 &\lesssim \Delta t^{-1} \Vert r \Vert_{\mathcal{Q}_T^n}^2 + \Delta t \Vert \partial_t r \Vert_{\mathcal{Q}_T^n}^2 \lesssim \Delta t \Vert u \Vert_{\mathcal{Q}_T^n}^2.
 \end{align}
 \end{subequations}
 Here, the constants depend on $\sigma$ through $\Vert \varphi \Vert_{W^{2,\infty}}$.
 Moreover, as the first estimate in \cref{eq_r_props_b} holds facet by facet, we may split $\mathcal{Q}_F \cap Q \subseteq \mathcal{Q}_F^{\text{gp}} \cup \mathcal{Q}_F^{\text{int}}$, apply \cref{lemma_gp_facet_jumps} on the stabilised facets and use $F \times I_n \subseteq \overline{Q^n}$ on the interior facets to obtain
 \begin{equation}
  \alpha \Vert h^{-\frac12} \fjump{r} \Vert_{\mathcal{Q}_F \cap Q}^2 \lesssim \Delta t^2 \big( \gamma_J^{-1} |u|_J^2 + \alpha \Vert h^{-\frac12} \fjump{u} \Vert_{\mathcal{Q}_F \cap Q}^2 \big) \lesssim \Delta t^2 \Vert u \Vert_{c,j}^2. \label{eq_r_jumps}
 \end{equation}
 We now estimate the five contributions to $(B_h+J)(u,r)$. For the volume convection term we obtain with Young's inequality (with weight $\Delta t^{\frac12}$) and \cref{lemma_gp_E_IpGP}
 \begin{align*}
  |C(u,r)| &\leq \Vert \w^\ast \cdo \nabla^\ast u \Vert_Q \Vert r \Vert_{\mathcal{Q}_T} \lesssim \Delta t^{\frac12} \left( \Delta t \Vert \w^\ast \cdo \nabla^\ast u \Vert_Q^2 + \Vert u \Vert_{\mathcal{Q}_T}^2 \right) \\
  &\lesssim \Delta t^{\frac12} \left( \Delta t \Vert \w^\ast \cdo \nabla^\ast u \Vert_Q^2 + \Vert u \Vert_{c,j}^2 \right),
 \end{align*}
 where $\Delta t \Vert \w^\ast \cdo \nabla^\ast u \Vert_Q^2 = \Delta t \Vert \mathcal{L} u \Vert_Q^2$ in Case C and $\Delta t \Vert \w^\ast \cdo \nabla^\ast u \Vert_Q^2 \lesssim \Delta t \Vert \mathcal{L} u \Vert_Q^2 + \Delta t \Vert \w \Vert_\infty^2 \Vert \nabla u \Vert_Q^2 \lesssim \tripplenorm{u}^2$ in Case D by \cref{eq:ass_caseD}.
 For the upwind terms we obtain with \cref{eq_r_props}
 \begin{align*}
  |U_h(u,r)| &\leq \tripplejump{u} \Big( \Vert |\w \cdo \n|^{\frac12} r_+ \Vert_{\mathcal{Q}_F \cap Q}^2 + \sum_{n=1}^N \Vert r_+ \Vert_{\Omega^{n-1}}^2 \Big)^{\frac12},
  \quad \sum_{n=1}^N \Vert r_+ \Vert_{\Omega^{n-1}}^2 \lesssim \Delta t \Vert u \Vert_{\mathcal{Q}_T}^2,
 \end{align*}
 and $\Vert |\w \cdo \n|^{\frac12} r_+ \Vert_{\mathcal{Q}_F \cap Q}^2 \lesssim \Vert \w \Vert_\infty h^{-1} \Delta t^2 \Vert u \Vert_{\mathcal{Q}_T}^2 \lesssim \Delta t \Vert u \Vert_{\mathcal{Q}_T}^2$ in Case C by the CFL-type condition. In Case D we estimate instead, using $\Vert \w \Vert_\infty^2 \Delta t \lesssim \alpha$,
 \begin{align*}
 (|\w \cdo \n| \fjump{u}, r_+)_{\mathcal{Q}_F \cap Q} &\leq \Vert \w \Vert_\infty \Vert \fjump{u} \Vert_{\mathcal{Q}_F \cap Q} \Vert r_+ \Vert_{\mathcal{Q}_F} \lesssim \Vert \w \Vert_\infty \Delta t \Vert h^{-\frac12} \fjump{u} \Vert_{\mathcal{Q}_F \cap Q} \Vert u \Vert_{\mathcal{Q}_T} \\ &\lesssim \Delta t^{\frac12} \big( \alpha \Vert h^{-\frac12} \fjump{u} \Vert_{\mathcal{Q}_F \cap Q}^2 + \Vert u \Vert_{\mathcal{Q}_T}^2 \big).
 \end{align*} In both cases we arrive at $|U_h(u,r)| \lesssim \Delta t^{\frac12} \Vert u \Vert_{c,j}^2$.
 For the diffusion terms we have with \cref{eq_r_props}, \cref{eq_r_jumps}, the trace inequality \cref{eq:std_trace} and \cref{eq_alpha_grad_full} (applied to $u$ and to $r$)
 \begin{align*}
  |A(u,r)| &\leq \alpha \Vert \nabla u \Vert_Q \Vert \nabla r \Vert_{\mathcal{Q}_T} \lesssim \Delta t \, \alpha \Vert \nabla u \Vert_{\mathcal{Q}_T}^2 \lesssim \Delta t \Vert u \Vert_{c,j}^2, \\
  |N_h(u,r)| &\leq \alpha \Vert h^{\frac12} \avg{\nabla u \cdo \n} \Vert_{\mathcal{Q}_F} \Vert h^{-\frac12} \fjump{r} \Vert_{\mathcal{Q}_F \cap Q} + \alpha \Vert h^{\frac12} \avg{\nabla r \cdo \n} \Vert_{\mathcal{Q}_F} \Vert h^{-\frac12} \fjump{u} \Vert_{\mathcal{Q}_F \cap Q} \\
  & \quad + \alpha \lambda \Vert h^{-\frac12} \fjump{u} \Vert_{\mathcal{Q}_F \cap Q} \Vert h^{-\frac12} \fjump{r} \Vert_{\mathcal{Q}_F \cap Q} \\
  & \lesssim \alpha^{\frac12} \Vert \nabla u \Vert_{\mathcal{Q}_T} \, \Delta t \Vert u \Vert_{c,j} + \Delta t \, \alpha^{\frac12} \Vert \nabla u \Vert_{\mathcal{Q}_T} \Vert u \Vert_c + \Vert u \Vert_c \, \Delta t \Vert u \Vert_{c,j} \lesssim \Delta t \Vert u \Vert_{c,j}^2, \\
  |J(u,r)| &\leq |u|_J |r|_J \lesssim \Delta t |u|_J^2.
 \end{align*}
 Collecting the estimates for the main part and the remainder yields the claim with $c_\sigma = \frac12 e^{-\sigma T} \min\{1, \sigma\}$.
\end{proof}

\subsection{Control of the derivative and inf-sup stability} \label{ssec:infsup}
To obtain control in the stronger norm $\tripplenorm{\cdot}$ we consider the test function $\Delta t \Lapp u \in W_h$, cf. \cref{eq:def_Lapp}. In Case D this is the (scaled) time derivative as in \cite{heimann2025discretizationerroranalysishigh}, in Case C the (scaled) discrete material derivative as in \cite{burman2025transport}.
\begin{lemma} \label{lemma_derivative_testing}
 There is a constant $C_{\ref{lemma_derivative_testing}} > 0$, independent of $h$, $\Delta t$ and $\alpha$, such that for all $u \in W_h$ (in both cases)
 \begin{equation}
  (B_h + J)(u, \Delta t \Lapp u) \geq \frac14 \Delta t \Vert \mathcal{L} u \Vert_Q^2 - C_{\ref{lemma_derivative_testing}} \Vert u \Vert_{c,j}^2. \label{eq_derivative_testing}
 \end{equation}
\end{lemma}
\begin{proof}
With $y := \Delta t \Lapp u$, we treat the five contributions $(C) := C(u,y)$, $(U) := U_h(u,y)$, $(A) := A(u,y)$, $(N) := N_h(u,y)$ and $(J) := J(u,y)$ of $(B_h + J)(u, y)$ separately, bounding each of them by its own constant $C_C$, $C_U$, $C_A$, $C_N$, $C_J > 0$, independent of $h$, $\Delta t$, $\alpha$ and $\epsilon$, where $\epsilon \in (0, 1]$ is chosen at the end. Throughout, we use $\Delta t \Vert \Lapp u \Vert_{\mathcal{Q}_T}^2 + |y|_J^2 \lesssim \Delta t \Vert \mathcal{L} u \Vert_Q^2 + \Vert u \Vert_{c,j}^2$ from \cref{eq_Lapp_tripplenorm} and $\alpha \Vert \nabla u \Vert_{\mathcal{Q}_T}^2 + \alpha \Vert \nabla y \Vert_{\mathcal{Q}_T}^2 \lesssim \Vert u \Vert_{c,j}^2$ from \cref{eq_alpha_grad_full}.

\underline{(C):}
By definition of $\mathcal{L}$ and $\Lapp$, Cauchy-Schwarz and Young's inequality
\begin{align*}
 (C) &= \Delta t ( \w^\ast \cdo \nabla^\ast u, \Lapp u )_Q = \Delta t \Vert \mathcal{L} u \Vert_Q^2 + \Delta t ( \mathcal{L} u, \Lapp u - \mathcal{L} u)_Q + \Delta t (\w^\ast \cdo \nabla^\ast u - \mathcal{L} u, \Lapp u)_Q.
\end{align*}
In Case C the last term vanishes and $\Delta t \Vert \Lapp u - \mathcal{L} u \Vert_Q^2 \lesssim \Delta t h^2 \Vert \nabla u \Vert_Q^2 \lesssim \Vert u \Vert_{c,j}^2$ by \cref{eq_Lapp_diff} and the proof of \cref{lemma_Lapp_tripplenorm}. In Case D the second term vanishes and $\Delta t \Vert \w^\ast \cdo \nabla^\ast u - \mathcal{L} u \Vert_Q^2 = \Delta t \Vert \w \cdo \nabla u \Vert_Q^2 \lesssim \alpha \Vert \nabla u \Vert_Q^2$ by \cref{eq:ass_caseD}. Hence, in both cases $(C) \geq \frac12 \Delta t \Vert \mathcal{L} u \Vert_Q^2 - C_C \Vert u \Vert_{c,j}^2$.

\underline{(U):}
By Cauchy-Schwarz and the definition of $\tripplejump{\cdot}$,
\begin{align*}
 |(U)| &\leq \tripplejump{u} \Big( \Delta t^2 \Vert |\w \cdo \n|^{\frac12} (\Lapp u)_+ \Vert_{\mathcal{Q}_F \cap Q}^2 + \Delta t^2 \sum_{n=1}^N \Vert (\Lapp u)_+ \Vert_{\Omega^{n-1}}^2 \Big)^{\frac12} \\
 &\leq \frac{1}{2 \epsilon} \tripplejump{u}^2 + \frac{\epsilon}{2} \Delta t^2 \Big( \Vert \w \Vert_\infty \Vert (\Lapp u)_+ \Vert_{\mathcal{Q}_F}^2 + \sum_{n=1}^N \Vert (\Lapp u)_+ \Vert_{\Omega^{n-1}}^2 \Big).
\end{align*}
In both cases, \cref{lemma_bound_tn_ag_temporal_vol} (first inequality) gives $\Delta t^2 \sum_n \Vert (\Lapp u)_+ \Vert_{\Omega^{n-1}}^2 \lesssim \Delta t \Vert \Lapp u \Vert_{\mathcal{Q}_T}^2$. In Case C, \cref{eq:std_trace} and the CFL-type condition moreover give $\Delta t^2 \Vert \w \Vert_\infty \Vert (\Lapp u)_+ \Vert_{\mathcal{Q}_F}^2 \lesssim \Vert \w \Vert_\infty \Delta t h^{-1} \cdot \Delta t \Vert \Lapp u \Vert_{\mathcal{Q}_T}^2 \lesssim \Delta t \Vert \Lapp u \Vert_{\mathcal{Q}_T}^2$. In Case D no CFL-type condition is available; here we apply the above splitting only to the contributions on $\Omega^{n-1}$ and bound the facet contribution directly, using $\Lapp u = \partial_t u$ and $\Vert \w \Vert_\infty^2 \Delta t \lesssim \alpha$,
\begin{align*}
 (|\w \cdo \n| \fjump{u}, \Delta t (\partial_t u)_+)_{\mathcal{Q}_F \cap Q} &\leq \Vert \w \Vert_\infty \Vert \fjump{u} \Vert_{\mathcal{Q}_F \cap Q} \Delta t \Vert (\partial_t u)_+ \Vert_{\mathcal{Q}_F} \\ &\lesssim \Vert \w \Vert_\infty \Delta t^{\frac12} \Vert h^{-\frac12} \fjump{u} \Vert_{\mathcal{Q}_F \cap Q} \cdot \Delta t^{\frac12} \Vert \partial_t u \Vert_{\mathcal{Q}_T} \\
 &\leq \frac{1}{2 \epsilon} \alpha \Vert h^{-\frac12} \fjump{u} \Vert_{\mathcal{Q}_F \cap Q}^2 + \frac{\epsilon}{2} \Delta t \Vert \partial_t u \Vert_{\mathcal{Q}_T}^2.
\end{align*} With \cref{lemma_Lapp_tripplenorm} we obtain in both cases $|(U)| \leq \frac{C_U}{\epsilon} \Vert u \Vert_{c}^2 + \epsilon C_U \tripplenorm{u}_j^2$.

\underline{(A):}
With Cauchy-Schwarz, Young's inequality and \cref{eq_alpha_grad_full},
\begin{equation*}
 |(A)| = |\alpha (\nabla u, \nabla y)_{Q}| \leq \frac{1}{2} \alpha \Vert \nabla u \Vert_{Q}^2 + \frac{1}{2} \alpha \Vert \nabla y \Vert_{\mathcal{Q}_T}^2 \leq C_A \Vert u \Vert_{c,j}^2.
\end{equation*}

\underline{(N):}
We apply Cauchy-Schwarz to each of the three contributions of $(N)$,
\begin{align*}
 [1] &:= ( \avg{ \alpha \nabla u \cdo \n}, \fjump{y})_{\mathcal{Q}_F \cap Q}, \qquad [2] := ( \avg{\alpha \nabla y \cdo \n}, \fjump{u })_{\mathcal{Q}_F \cap Q}, \\
 [3] &:= ( \alpha \lambda h^{-1} \fjump{u}, \fjump{y})_{\mathcal{Q}_F \cap Q}.
\end{align*} The arising factors are bounded as follows: By \cref{eq:std_trace} and \cref{eq_alpha_grad_full}, $\alpha \Vert h^{\frac12} \avg{\nabla u \cdo \n} \Vert_{\mathcal{Q}_F}^2 \lesssim \alpha \Vert \nabla u \Vert_{\mathcal{Q}_T}^2 \lesssim \Vert u \Vert_{c,j}^2$ and similarly $\alpha \Vert h^{\frac12} \avg{\nabla y \cdo \n} \Vert_{\mathcal{Q}_F}^2 \lesssim \alpha \Vert \nabla y \Vert_{\mathcal{Q}_T}^2 \lesssim \Vert u \Vert_{c,j}^2$. By \cref{lemma_y_jumps}, $\alpha \Vert h^{-\frac12} \fjump{y} \Vert_{\mathcal{Q}_F \cap Q}^2 \lesssim \tripplenorm{u}_j^2$. Finally, $\alpha \Vert h^{-\frac12} \fjump{u} \Vert_{\mathcal{Q}_F \cap Q}^2 \leq \Vert u \Vert_c^2$. Hence, $|[1]| + |[3]| \leq \frac{C}{\epsilon} \Vert u \Vert_{c,j}^2 + \epsilon C \tripplenorm{u}_j^2$ and $|[2]| \leq C \Vert u \Vert_{c,j}^2$, so that, using $\epsilon \leq 1$ for $[2]$, $|(N)| \leq \frac{C_N}{\epsilon} \Vert u \Vert_{c,j}^2 + \epsilon C_N \tripplenorm{u}_j^2$.

\underline{(J):}
By Cauchy-Schwarz, Young's inequality and \cref{lemma_Lapp_tripplenorm}, $|(J)| \leq |u|_J |y|_J \leq \frac{C_J}{\epsilon} |u|_J^2 + \epsilon C_J \tripplenorm{u}_j^2$.

\underline{Putting it all together:}
Summing the five bounds, using $\Vert u \Vert_c, |u|_J \leq \Vert u \Vert_{c,j}$ and $\tripplenorm{u}_j^2 = \Delta t \Vert \mathcal{L} u \Vert_Q^2 + \Vert u \Vert_{c,j}^2$, cf. \cref{eq_def_tripplenorm}, we obtain
\begin{align*}
(B_h + J)(u, y) \geq{}& \frac12 \Delta t \Vert \mathcal{L} u \Vert_Q^2 - C_C \Vert u \Vert_{c,j}^2 && \leftarrow (C) \\
 & - \frac{C_U}{\epsilon} \Vert u \Vert_{c,j}^2 - \epsilon C_U \tripplenorm{u}_j^2 && \leftarrow (U) \\
 & - C_A \Vert u \Vert_{c,j}^2 && \leftarrow (A) \\
 & - \frac{C_N}{\epsilon} \Vert u \Vert_{c,j}^2 - \epsilon C_N \tripplenorm{u}_j^2 && \leftarrow (N) \\
 & - \frac{C_J}{\epsilon} \Vert u \Vert_{c,j}^2 - \epsilon C_J \tripplenorm{u}_j^2 && \leftarrow (J) \\
 ={}& \Big( \frac12 - \epsilon C_\ast \Big) \Delta t \Vert \mathcal{L} u \Vert_Q^2
   - \Big( C_C + C_A + \big( \epsilon + \tfrac{1}{\epsilon} \big) C_\ast \Big) \Vert u \Vert_{c,j}^2
\end{align*}
with $C_\ast := C_U + C_N + C_J$, where we may assume $C_\ast \geq \frac14$ (enlarging $C_U$ if necessary). Choosing $\epsilon = (4 C_\ast)^{-1} \leq 1$ yields the claim with $C_{\ref{lemma_derivative_testing}} = C_C + C_A + 4 C_\ast^2 + \frac14$.
\end{proof}

In the stability proof, we will exploit the compatibility of the $\tripplenorm{\cdot}_j$ norms between a function $u$ and the specially tailored test functions.
\begin{lemma} \label{lemma_test_fun_bnd}
 For all $u \in W_h$ it holds
 \begin{equation}
  \tripplenorm{\Delta t \Lapp u}_j \lesssim \tripplenorm{u}_j \qquad \text{ and } \qquad \tripplenorm{\Pi_t(e^{-\sigma t} u)}_j \lesssim \tripplenorm{u}_j, \label{eq_test_fun_bnd}
 \end{equation}
 where the constant in the second estimate depends on $\sigma$.
\end{lemma}
\begin{proof}
 We bound the summands of the squared norms individually.
 For $y := \Delta t \Lapp u$, the inverse inequalities \cref{eq:std_inverse} and, in Case C, the CFL-type condition yield
 \begin{equation*}
  \Delta t \Vert \mathcal{L} y \Vert_Q^2 \lesssim \Delta t^3 \left( \Delta t^{-2} + \Vert \w \Vert_\infty^2 h^{-2} \right) \Vert \Lapp u \Vert_{\mathcal{Q}_T}^2 \lesssim \Delta t \Vert \Lapp u \Vert_{\mathcal{Q}_T}^2,
 \end{equation*}
 where the term with $\w$ is absent in Case D, and $\Vert y \Vert_Q^2 \leq \Delta t \cdot \Delta t \Vert \Lapp u \Vert_{\mathcal{Q}_T}^2$. The temporal traces in $\tripplejump{y}^2$ are bounded by $\Delta t \Vert \Lapp u \Vert_{\mathcal{Q}_T}^2$ by \cref{lemma_bound_tn_ag_temporal_vol}, the spatial jumps in $\tripplejump{y}^2$ and in $\alpha |y|_{1,h}^2$ by \cref{lemma_y_jumps}, and $\Delta t \Vert \Lapp u \Vert_{\mathcal{Q}_T}^2$, $\alpha \Vert \nabla y \Vert_Q^2$ and $|y|_J^2$ by \cref{lemma_Lapp_tripplenorm}.

 For $u_w := \Pi_t(e^{-\sigma t} u) = \varphi u + r$ with, as in the proof of \cref{lemma_weighted_symm_testing}, $\varphi := e^{-\sigma t}$ and $r := (\Pi_t - \mathrm{Id})(\varphi u)$, the norms of $\varphi u$ are bounded by those of $u$, where $\mathcal{L} (\varphi u) = \varphi (\mathcal{L} u - \sigma u)$. For $r$, the estimates \cref{eq_r_props}, \cref{eq_r_jumps} and \cref{eq_alpha_grad_full} bound the $L^2$, trace, diffusion and Ghost penalty parts by $\Delta t \Vert u \Vert_{\mathcal{Q}_T}^2 + \Delta t^2 \Vert u \Vert_{c,j}^2$. The spatial part of $\tripplejump{r}^2$ is bounded by \cref{lemma_w_jumps}: in Case D by $\alpha \Vert h^{-\frac12} \fjump{r} \Vert_{\mathcal{Q}_F \cap Q}^2$ and in Case C by $\Delta t^{-1} \Vert r \Vert_{\mathcal{Q}_T}^2 \lesssim \Delta t \Vert u \Vert_{\mathcal{Q}_T}^2$. For the derivative term, $\Delta t \Vert \partial_t r \Vert_Q^2 \lesssim \Delta t \Vert u \Vert_{\mathcal{Q}_T}^2$, and in Case C additionally $\Delta t \Vert \w \Vert_\infty^2 \Vert \nabla r \Vert_Q^2 \lesssim \Delta t^3 \Vert \w \Vert_\infty^2 h^{-2} \Vert u \Vert_{\mathcal{Q}_T}^2 \lesssim \Delta t \Vert u \Vert_{\mathcal{Q}_T}^2$ by the CFL-type condition. Finally, $\Delta t \Vert u \Vert_{\mathcal{Q}_T}^2 \lesssim \Vert u \Vert_{c,j}^2$ by \cref{lemma_gp_E_IpGP} and \cref{eq:ass_gp_scaling}.
\end{proof}

We are now in the position to prove the main stability result.
\begin{theorem}[Inf-sup stability] \label{thm_stability}
 Let the assumptions of \cref{ssec:prelim} hold (Case D or Case C), $\lambda \geq \lambda_0$, and let $\Delta t$ be sufficiently small, i.e. $\Delta t \leq \Delta t_0$ with $\Delta t_0$ depending on $\sigma$, $T$ and the constants of \cref{lemma_weighted_symm_testing,lemma_derivative_testing}, but not on $h$ or $\alpha$. Then, the bilinear form $B_h + J$ is inf-sup stable on $W_h$ w.r.t. the norm $\tripplenorm{\cdot}_j$, i.e. there is $C_{\tilde v} > 0$ such that for all $u \in W_h$ there exists $\tilde v(u) \in W_h$ with
 \begin{equation}
  \tripplenorm{\tilde v(u)}_j \leq C_{\tilde v} \tripplenorm{u}_j \qquad \text{ and } \qquad (B_h + J)(u, \tilde v(u)) \geq \frac18 \tripplenorm{u}_j^2. \label{eq_inf_sup_two}
 \end{equation}
 In particular, $\tilde v(u) \neq 0$ for $u \neq 0$, and with $C_{\text{stab}} := (8 C_{\tilde v})^{-1}$ we have $(B_h + J)(u, \tilde v(u)) \geq C_{\text{stab}} \tripplenorm{u}_j \tripplenorm{\tilde v(u)}_j$, i.e.
 \begin{equation}
  \inf_{u \in W_h \setminus \{0\}} \, \sup_{v \in W_h \setminus \{0\}} \frac{(B_h + J)(u, v)}{\tripplenorm{u}_j \tripplenorm{v}_j} \geq C_{\text{stab}}. \label{eq_inf_sup}
 \end{equation}
 $C_{\text{stab}}$ depends on $\sigma$ and $T$ through $c_\sigma \simeq e^{-\sigma T} \min\{1,\sigma\}$, but is independent of $h$, $\Delta t$ and $\alpha$.
 In particular, the discrete problem \cref{eq_discrete_problem} has a unique solution $u_h \in W_h$.
\end{theorem}
\begin{proof}
 Fix $u \in W_h$ and set $\tilde v(u) := K \Pi_t(e^{-\sigma t} u) + \Delta t \Lapp u \in W_h$ with $K := (C_{\ref{lemma_derivative_testing}} + \frac14) / c_\sigma$. By \cref{lemma_weighted_symm_testing} and \cref{lemma_derivative_testing}
 \begin{align*}
  (B_h + J)(u, \tilde v(u)) &\geq K c_\sigma \Vert u \Vert_{c,j}^2 - K C \Delta t^{\frac12} \tripplenorm{u}_j^2 + \frac14 \Delta t \Vert \mathcal{L} u \Vert_Q^2 - C_{\ref{lemma_derivative_testing}} \Vert u \Vert_{c,j}^2 \\
  &= \frac14 \left( \Vert u \Vert_{c,j}^2 + \Delta t \Vert \mathcal{L} u \Vert_Q^2 \right) - K C \Delta t^{\frac12} \tripplenorm{u}_j^2
  \geq \frac18 \tripplenorm{u}_j^2
 \end{align*}
 for $\Delta t \leq \Delta t_0 := \min\{ \sigma^{-1}, (8 K C)^{-2} \}$, where $\sigma \Delta t \leq 1$ is needed for \cref{lemma_weighted_symm_testing} and the second restriction ensures $K C \Delta t^{\frac12} \leq \frac18$. By \cref{lemma_test_fun_bnd}, $\tripplenorm{\tilde v(u)}_j \leq (K C' + C'') \tripplenorm{u}_j$, which is the first estimate in \cref{eq_inf_sup_two} with $C_{\tilde v} = K C' + C''$. Now recall that $\tripplenorm{\cdot}_j$ is a norm on $W_h$ (note that $(\Vert \cdot \Vert_Q^2 + |\cdot|_J^2)^{\frac12}$ is a norm on $W_h$ by \cref{lemma_gp_E_IpGP}). Hence, for $u \neq 0$ the second estimate in \cref{eq_inf_sup_two} gives $(B_h + J)(u, \tilde v(u)) > 0$, so that $\tilde v(u) \neq 0$, and combining both estimates yields \cref{eq_inf_sup}. Unique solvability follows as $W_h$ is finite dimensional.
\end{proof}
\begin{remark}[On the role of the $L^2$ norm and the weighting]
 In Case D, the weighting with $e^{-\sigma t}$ and the $L^2(Q)$ contribution in $\Vert \cdot \Vert_c$ could be avoided by using the supremiser $\tilde v(u) = K u + \Delta t \partial_t u$; $L^2$ control then has to be recovered separately, cf. \cite{heimann2025discretizationerroranalysishigh}. In Case C, the terms stemming from the approximation $\wapp$ of $\w$ cannot be absorbed without $L^2$ control, so that the weighting is essential, cf. \cite{burman2025transport}. For a unified presentation we use the weighted test function in both cases.
\end{remark}

\subsection{Continuity} \label{ssec:continuity}
\begin{lemma} \label{lemma_continuity}
 For all $u \in W_h + H^2(Q)$ and $v \in W_h$ it holds
 \begin{equation}
  B_h(u,v) \lesssim \tripplenorm{u}_\ast \tripplenorm{v}_j \quad \text{ and } \quad J(u,v) \leq |u|_J |v|_J \text{ for } u,v \in W_h. \label{eq_continuity}
 \end{equation}
\end{lemma}
\begin{proof}
 The second claim is Cauchy-Schwarz. For the first claim we write the convective part in the form of \cref{lemma_B_Bmc} and apply Cauchy-Schwarz to each term:
 \begin{align*}
  |(C + U_h)(u,v)| &\leq |(u, \w^\ast \cdo \nabla^\ast v)_Q| + \Vert |\w \cdo \n|^{\frac12} u_- \Vert_{\mathcal{Q}_F \cap Q} \Vert |\w \cdo \n|^{\frac12} \fjump{v} \Vert_{\mathcal{Q}_F \cap Q} \\
  & \quad + \sum_{n=1}^{N-1} \Vert u_- \Vert_{\Omega^n} \Vert \fjump{v} \Vert_{\Omega^n} + \Vert u_- \Vert_{\Omega^N} \Vert v_- \Vert_{\Omega^N}.
 \end{align*}
 By Cauchy-Schwarz in $n$, the temporal terms are bounded by $\big( \sum_{n=1}^N \Vert u_- \Vert_{\Omega^n}^2 \big)^{\frac12} \tripplejump{v}$.
 In Case C, $(u, \w^\ast \cdo \nabla^\ast v)_Q \leq \Delta t^{-\frac12} \Vert u \Vert_Q \, \Delta t^{\frac12} \Vert \mathcal{L} v \Vert_Q$ and $\Vert \w \Vert_\infty \lesssim h / \Delta t$, so that $\Vert |\w \cdo \n|^{\frac12} u_- \Vert_{\mathcal{Q}_F \cap Q}^2 \lesssim \frac{h}{\Delta t} \Vert u_- \Vert_{\mathcal{Q}_F \cap Q}^2$. In Case D, we use $\Vert \w \Vert_\infty \alpha^{-\frac12} \lesssim \Delta t^{-\frac12}$ and $\Vert \w \Vert_\infty^2 h \alpha^{-1} \lesssim h / \Delta t$, cf. \cref{eq:ass_caseD}, in
 \begin{align*}
  (u, \w^\ast \cdo \nabla^\ast v)_Q &\leq \Delta t^{-\frac12} \Vert u \Vert_Q \, \Delta t^{\frac12} \Vert \partial_t v \Vert_Q + \Vert \w \Vert_\infty \alpha^{-\frac12} \Vert u \Vert_Q \, \alpha^{\frac12} \Vert \nabla v \Vert_Q, \\
  (|\w \cdo \n| u_-, \fjump{v})_{\mathcal{Q}_F \cap Q} &\leq \Vert \w \Vert_\infty h^{\frac12} \alpha^{-\frac12} \Vert u_- \Vert_{\mathcal{Q}_F \cap Q} \, \alpha^{\frac12} \Vert h^{-\frac12} \fjump{v} \Vert_{\mathcal{Q}_F \cap Q}.
 \end{align*}
 In both cases all factors involving $u$ are bounded by $\tripplenorm{u}_\ast$ and all factors involving $v$ by $\tripplenorm{v}$.
 For the diffusive part we have
 \begin{align*}
  |(A + N_h)(u,v)| \leq \ & \alpha \Vert \nabla u \Vert_Q \Vert \nabla v \Vert_Q + \alpha \Vert h^{\frac12} \avg{\nabla u \cdo \n} \Vert_{\mathcal{Q}_F \cap Q} \Vert h^{-\frac12} \fjump{v} \Vert_{\mathcal{Q}_F \cap Q} \\
  & + \alpha \Vert h^{\frac12} \avg{\nabla v \cdo \n} \Vert_{\mathcal{Q}_F} \Vert h^{-\frac12} \fjump{u} \Vert_{\mathcal{Q}_F \cap Q} \\
  & + \alpha \lambda \Vert h^{-\frac12} \fjump{u} \Vert_{\mathcal{Q}_F \cap Q} \Vert h^{-\frac12} \fjump{v} \Vert_{\mathcal{Q}_F \cap Q},
 \end{align*}
 where $\alpha \Vert h^{\frac12} \avg{\nabla v \cdo \n} \Vert_{\mathcal{Q}_F}^2 \lesssim \alpha \Vert \nabla v \Vert_{\mathcal{Q}_T}^2 \lesssim \alpha \Vert \nabla v \Vert_Q^2 + \gamma_J^{-1} |v|_J^2 \lesssim \tripplenorm{v}_j^2$ by \cref{eq:std_trace} and \cref{eq_alpha_grad_full}. This is the only place where the Ghost penalty part of $\tripplenorm{v}_j$ is needed.
\end{proof}

\section{A priori error analysis} \label{sect_a_pr_est}
\subsection{Consistency}
\begin{lemma} \label{lemma_consistency}
 Let $u \in H^2(Q)$ be the solution of \cref{strongformproblem}. Then, for all $v \in W_h$ it holds
 \begin{equation}
  B_h(u,v) = f_h(v). \label{eq_consistency}
 \end{equation}
\end{lemma}
\begin{proof}
 As $u$ is continuous in time and across spatial facets, all jump terms of $u$ vanish, except for the term on $\Omega^0$ where $\fjump{u} = u_+ = u_0$ by convention. Hence, $U_h(u,v) = (u_0, v_+)_{\Omega^0}$, and $N_h(u,v) = (\alpha \nabla u \cdo \n, \fjump{v})_{\mathcal{Q}_F \cap Q}$. Elementwise partial integration in space on $Q^n$ yields, with the outer normal $\n_T$ of the element $T$ and $\n_{\partial \Omega}$ the outer normal of $\Omega(t)$,
 \begin{align*}
  A(u,v) &= - \alpha (\Delta u, v)_Q + \alpha \sum_{n} \sum_{Q_T^n \in \mathcal{Q}_T^n} (\nabla u \cdo \n_T, v)_{\partial_s (Q_T^n \cap Q^n)} \\
  &= - \alpha (\Delta u, v)_Q + \alpha (\nabla u \cdo \n_{\partial \Omega}, v)_{\partial_s Q} - (\alpha \nabla u \cdo \n, \fjump{v})_{\mathcal{Q}_F \cap Q},
 \end{align*}
 where $\partial_s Q := \bigcup_{t \in (0,T)} \partial \Omega(t) \times \{t\}$ denotes the lateral boundary of $Q$, as for each interior facet the upwind element has outer normal $\n$ (with trace $v_-$) and the downwind element has outer normal $-\n$ (with trace $v_+$), so that $v_- \nabla u \cdo \n - v_+ \nabla u \cdo \n = - \nabla u \cdo \n \fjump{v}$. The boundary term vanishes due to \cref{strongformproblem:2} (and is absent for $\alpha = 0$). Hence,
 $B_h(u,v) = (\partial_t u + \w \cdo \nabla u - \alpha \Delta u, v)_Q + (u_0, v_+)_{\Omega^0} = (f,v)_Q + (u_0, v_+)_{\Omega^0} = f_h(v)$ by \cref{strongformproblem:1} and \cref{strongformproblem:3}.
\end{proof}

\subsection{Abstract error bound}
Combining stability, continuity and consistency in the usual way for nonconforming methods, cf. e.g. \cite[Thm. 1.35]{di2011mathematical}, we obtain the following Strang-type result.
\begin{proposition} \label{prop_strang}
 Let the assumptions of \cref{thm_stability} hold, let $u \in H^2(Q)$ be the solution of \cref{strongformproblem} and $u_h \in W_h$ the solution of \cref{eq_discrete_problem}. Then,
 \begin{equation}
  \tripplenorm{u - u_h} \lesssim \inf_{w_h \in W_h} \left( \tripplenorm{u - w_h}_\ast + |w_h|_J \right). \label{eq_strang}
 \end{equation}
\end{proposition}
\begin{proof}
 Fix $w_h \in W_h$ and set $e_h := u_h - w_h \in W_h$. By \cref{lemma_consistency} and \cref{eq_discrete_problem}, $(B_h + J)(e_h, v) = B_h(u - w_h, v) - J(w_h, v)$ for all $v \in W_h$. Hence, by \cref{thm_stability} and \cref{lemma_continuity}
 \begin{align*}
  C_{\text{stab}} \tripplenorm{e_h}_j &\leq \sup_{v \in W_h \setminus \{0\}} \frac{(B_h + J)(e_h, v)}{\tripplenorm{v}_j} = \sup_{v \in W_h \setminus \{0\}} \frac{B_h(u - w_h, v) - J(w_h, v)}{\tripplenorm{v}_j} \\ &\lesssim \tripplenorm{u - w_h}_\ast + |w_h|_J.
 \end{align*}
 The triangle inequality $\tripplenorm{u - u_h} \leq \tripplenorm{u - w_h} + \tripplenorm{e_h}_j$ and $\tripplenorm{\cdot} \leq \tripplenorm{\cdot}_\ast$ yield the claim.
\end{proof}

\subsection{Interpolation}
It remains to bound the right hand side of \cref{eq_strang} for a suitable choice $w_h = I_h u$. As the active mesh region $Q_{\mathcal{E}} \supset Q$ is larger than the physical domain, we need an extension of the solution. We use the degree-independent Sobolev extension operator of \cite{S70} (in space-time), and assume the following regularity.
\begin{assumption} \label[assumption]{ass:regularity}
 The solution $u$ of \cref{strongformproblem} admits an extension $u^e \in H^{m}(\tilde \Omega \times (0,T))$, $u^e|_Q = u$, with $\Vert u^e \Vert_{H^m(\tilde \Omega \times (0,T))} \lesssim \Vert u \Vert_{H^m(Q)}$ for the regularity index $m = k_{\max} + 2$, $k_{\max} := \max\{k_s, k_t\}$, cf. \cref{thm_apriori} below.
\end{assumption}
To distinguish spatial and temporal regularity one could use $t$-anisotropic Sobolev spaces as in \cite{heimann2025discretizationerroranalysishigh, burman_preuss_25}. For ease of presentation we assume maximal regularity when compared to the discrete orders.

As interpolation operator we choose $I_h := I_s \otimes I_t$, the tensor product of a spatial and a temporal interpolation, where $I_t$ denotes the interpolation into $\mathcal{P}^{k_t}(I_n)$ on each time slab (acting pointwise in space, e.g. the $L^2(I_n)$-projection $\Pi_t$) and $I_s$ the elementwise $L^2(T)$-projection into $\mathcal{P}^{k_s}(T)$ (acting pointwise in time), restricted to the active elements. We assume $k_s \geq 1$ (for $k_s = 0$ the estimates involving the spatial Hessian below have to be modified). By construction $I_h u^e \in W_h$, and the following estimates are standard, cf. \cite[Thm. A5, Lemma A7]{burman_preuss_25}, \cite[Sect. 6.2]{heimann2025discretizationerroranalysishigh}: For $u^e$ as in \cref{ass:regularity}, and with $\Vert \cdot \Vert_{H^m} := \Vert \cdot \Vert_{H^m(\tilde \Omega \times (0,T))}$,
\begin{subequations} \label{eq_interpol}
\begin{align}
  \| u^e - I_h u^e \|_{\mathcal{Q}_T} &\lesssim \Delta t^{k_t + 1} \| u^e \|_{H^{k_t+1}} + h^{k_s +1} \| u^e \|_{H^{k_s+1}}, \label{eq_interpol_L2} \\
  \| \partial_t (u^e - I_h u^e) \|_{\mathcal{Q}_T} &\lesssim \Delta t^{k_t} \| u^e \|_{H^{k_t+1}} + h^{k_s +1} \| u^e \|_{H^{k_s+2}}, \label{eq_interpol_dt} \\
  \| \nabla (u^e - I_h u^e) \|_{\mathcal{Q}_T} + h \| D^2 (u^e - I_h u^e) \|_{\mathcal{Q}_T} &\lesssim \Delta t^{k_t+1} \| u^e \|_{H^{k_t+2}} + h^{k_s} \| u^e \|_{H^{k_s+1}}, \label{eq_interpol_grad}
\end{align}
and, for the spatial facet jumps and the Ghost penalty jumps, exploiting that $I_t u^e$ is continuous in space,
\begin{align}
 \Vert h^{-\frac12} \fjump{I_h u^e} \Vert_{\mathcal{Q}_F}^2 + h^{-2} \sum_{n=1}^N j_h^n(I_h u^e, I_h u^e) &\lesssim h^{2 k_s} \Vert u^e \Vert_{H^{k_s+1}}^2. \label{eq_interpol_jumps}
\end{align}
\end{subequations}
Here $D^2$ denotes the spatial Hessian. For the second order term in \cref{eq_interpol_grad} we use the splitting $u^e - I_h u^e = (u^e - I_s u^e) + I_s (u^e - I_t u^e)$: the first part is a standard spatial interpolation error, and for the second part the inverse inequality, the local $H^1$-stability of the $L^2(T)$-projection on shape regular elements and $\nabla I_t = I_t \nabla$ give $h \Vert D^2 I_s (u^e - I_t u^e) \Vert_{Q_T^n} \lesssim \Vert \nabla I_s (u^e - I_t u^e) \Vert_{Q_T^n} \lesssim \Vert (\mathrm{Id} - I_t) \nabla u^e \Vert_{Q_T^n} \lesssim \Delta t^{k_t+1} \Vert u^e \Vert_{H^{k_t+2}(Q_T^n)}$.
For \cref{eq_interpol_jumps} we refer to the proof of \cite[Prop. 3.26]{preuss18}: on each facet patch $\omega_F$, the two polynomials $I_h u^e|_{T_i}$, $i=1,2$, and their extensions $\mathcal{E}^p$ are both approximations of $I_t u^e$ of order $h^{k_s+1}$ in $L^2(\omega_F \times I_n)$, so that $\Vert \jump{I_h u^e}_{\omega_F} \Vert_{\omega_F \times I_n} \lesssim h^{k_s+1} \Vert I_t u^e \Vert_{L^2(I_n; H^{k_s+1}(\omega_F))} \lesssim h^{k_s+1} \Vert u^e \Vert_{L^2(I_n; H^{k_s+1}(\omega_F))}$ by the $L^2$-stability of $I_t$ (applied to the spatial derivatives). For the facet jumps we write $\fjump{I_h u^e}|_{F \times I_n} = \pm \jump{I_h u^e}_{\omega_F}|_{T_1}|_{F \times I_n}$ as in the proof of \cref{lemma_gp_facet_jumps} and apply the trace inequality \cref{eq:std_trace} to this polynomial.

\begin{lemma} \label{lemma_approx}
 Let $u$ satisfy \cref{ass:regularity} and $k := \min\{k_s, k_t\}$. Then, under the assumptions of Case D,
 \begin{align}
  \tripplenorm{u - I_h u^e}_\ast + |I_h u^e|_J &\lesssim \left( \alpha^{\frac12} \big( h^{k_s} + \Delta t^{k_t + 1} \big) + \Delta t^{k_t + \frac12} \right) \Vert u \Vert_{H^{k_{\max}+2}(Q)}, \label{eq_approx_caseD}
 \end{align}
 and under the assumptions of Case C,
 \begin{align}
  \tripplenorm{u - I_h u^e}_\ast + |I_h u^e|_J &\lesssim \left( h^{k + \frac12} + \alpha^{\frac12} \big( h^{k_s} + h^{k_t + 1} \big) \right) \Vert u \Vert_{H^{k_{\max}+2}(Q)}. \label{eq_approx_caseC}
 \end{align}
\end{lemma}
\begin{proof}
 Let $e := u^e - I_h u^e$ and note that $\tripplenorm{u - I_h u^e}_\ast$ only involves $e$ on $Q$ and the physical parts of the facets, so that we may bound all terms by the respective norms on $\mathcal{Q}_T$. We abbreviate $\Vert u^e \Vert := \Vert u^e \Vert_{H^{k_{\max}+2}} \lesssim \Vert u \Vert_{H^{k_{\max}+2}(Q)}$ and go through the summands of $\tripplenorm{e}_\ast^2$, cf. \cref{eq:norms}, \cref{eq_def_tripplenorm_ast}.

 \underline{Volume terms:} By \cref{eq_interpol_L2}, $(1 + \Delta t^{-1}) \Vert e \Vert_Q^2 \lesssim \Delta t^{-1} ( \Delta t^{2 k_t + 2} + h^{2 k_s + 2} ) \Vert u^e \Vert^2$. By \cref{eq_interpol_dt}, \cref{eq_interpol_grad},
 \begin{align*}
 \Delta t \Vert \mathcal{L} e \Vert_Q^2 &\lesssim \Delta t \Vert \partial_t e \Vert_{\mathcal{Q}_T}^2 + \Delta t \Vert \w \Vert_\infty^2 \Vert \nabla e \Vert_{\mathcal{Q}_T}^2 \\ &\lesssim \left( \Delta t^{2 k_t + 1} + \Delta t h^{2 k_s + 2} + \Delta t \Vert \w \Vert_\infty^2 h^{2 k_s} + \Vert \w \Vert_\infty^2 \Delta t^{2k_t+3} \right) \Vert u^e \Vert^2,
 \end{align*}
 where the last two terms only appear in Case C and are bounded by $h^{2k_s+1}$ and $h^{2k_t+3}$ there, by $h \simeq \Vert \w \Vert_\infty \Delta t$. By \cref{eq_interpol_grad} and \cref{eq_interpol_jumps}, $\alpha |e|_{1,h}^2 \lesssim \alpha ( h^{2 k_s} + \Delta t^{2 k_t + 2} ) \Vert u^e \Vert^2$.

 \underline{Trace terms:} By \cref{eq:cont_trace}, \cref{eq_interpol_L2} and \cref{eq_interpol_dt}, we have
 $\sum_{n} ( \Vert e_+ \Vert_{\Omega^{n-1}}^2 + \Vert e_- \Vert_{\Omega^n}^2 ) \lesssim \Delta t^{-1} \Vert e \Vert_{\mathcal{Q}_T}^2 + \Delta t \Vert \partial_t e \Vert_{\mathcal{Q}_T}^2 \lesssim ( \Delta t^{2k_t+1} + \Delta t^{-1} h^{2 k_s+2} ) \Vert u^e \Vert^2$, which bounds the temporal parts of $\tripplejump{e}^2$ and $\sum_n \Vert e_- \Vert_{\Omega^n}^2$. For the spatial facet terms we have by \cref{eq_interpol_jumps} and \cref{eq:cont_trace}
 \begin{align*}
  \Vert |\w \cdo \n|^{\frac12} \fjump{e} \Vert_{\mathcal{Q}_F \cap Q}^2 &\lesssim \Vert \w \Vert_\infty h^{2 k_s + 1} \Vert u^e \Vert^2, \\
  \frac{h}{\Delta t} \Vert e_\pm \Vert_{\mathcal{Q}_F}^2 &\lesssim \Delta t^{-1} \left( \Vert e \Vert_{\mathcal{Q}_T}^2 + h^2 \Vert \nabla e \Vert_{\mathcal{Q}_T}^2 \right) \lesssim \Delta t^{-1} \left( \Delta t^{2 k_t + 2} + h^{2 k_s + 2} \right) \Vert u^e \Vert^2, \\
  \alpha \Vert h^{\frac12} \avg{\nabla e \cdo \n} \Vert_{\mathcal{Q}_F}^2 &\lesssim \alpha \left( \Vert \nabla e \Vert_{\mathcal{Q}_T}^2 + h^2 \Vert D^2 e \Vert_{\mathcal{Q}_T}^2 \right) \lesssim \alpha \left( h^{2 k_s} + \Delta t^{2 k_t + 2} \right) \Vert u^e \Vert^2.
 \end{align*}

 \underline{Ghost penalty term:} By \cref{eq_interpol_jumps} and the definition of $c_J$,
 $$|I_h u^e|_J^2 = c_J \sum_n j_h^n(I_h u^e, I_h u^e) \lesssim \gamma_J ( \alpha + \Vert \w \Vert_\infty h ) h^{2 k_s} \Vert u^e \Vert^2.$$

 \underline{Collecting the terms:} In Case D we use $\Vert \w \Vert_\infty h \lesssim \alpha$, $h^2 \lesssim \alpha \Delta t$ and $\Delta t \leq 1$ (so that $\Delta t^{-1} h^{2k_s+2} \lesssim \alpha h^{2 k_s}$ and $\Delta t h^{2k_s+2} \leq h^{2k_s+2} \lesssim \alpha h^{2 k_s}$), which yields $\tripplenorm{e}_\ast^2 + |I_h u^e|_J^2 \lesssim ( \alpha h^{2 k_s} + \alpha \Delta t^{2 k_t + 2} + \Delta t^{2 k_t + 1} ) \Vert u^e \Vert^2$. In Case C we use $h \simeq \Vert \w \Vert_\infty \Delta t$ together with $\Vert \w \Vert_\infty \simeq 1$ (so that $\Delta t \simeq h$) and $h \lesssim 1$, cf. \cref{ssec:prelim}, so that all terms are bounded by $( h^{2 k + 1} + \alpha h^{2 k_s} + \alpha h^{2 k_t + 2} ) \Vert u^e \Vert^2$. Note that the diffusion-weighted temporal error $\alpha \Delta t^{2 k_t + 2}$ cannot be dropped in general: it is dominated by $\Delta t^{2k_t+1}$ only if $\alpha \Delta t \lesssim 1$, and (in Case C) by $h^{2k+1}$ if $\alpha \lesssim \Vert \w \Vert_\infty h$.
\end{proof}

\subsection{A priori error bound}
\begin{theorem}[A priori error bound] \label{thm_apriori}
 Let the assumptions of \cref{thm_stability} and \cref{ass:regularity} hold, $k_s \geq 1$, let $u$ be the solution of \cref{strongformproblem} and $u_h \in W_h$ the solution of \cref{eq_discrete_problem}. Then, in Case D (diffusion dominated, \cref{eq:ass_caseD}),
 \begin{equation}
  \tripplenorm{u - u_h} \lesssim \left( \alpha^{\frac12} \big( h^{k_s} + \Delta t^{k_t + 1} \big) + \Delta t^{k_t + \frac12} \right) \Vert u \Vert_{H^{k_{\max}+2}(Q)}, \label{eq_apriori_caseD}
 \end{equation}
 and in Case C (convection dominated, \cref{eq:ass_caseC}), with $k = \min\{k_s,k_t\}$,
 \begin{equation}
  \tripplenorm{u - u_h} \lesssim \big( h^{k + \frac12} + \alpha^{\frac12} ( h^{k_s} + h^{k_t + 1} ) \big) \Vert u \Vert_{H^{k_{\max}+2}(Q)}. \label{eq_apriori_caseC}
 \end{equation}
 In particular, if $\alpha \lesssim \Vert \w \Vert_\infty h$ in Case C, we have $\tripplenorm{u - u_h} \lesssim h^{k + \frac12} \Vert u \Vert_{H^{k_{\max}+2}(Q)}$, and if $\alpha \Delta t \lesssim 1$ (e.g. $\alpha \lesssim 1$) in Case D, we have $\tripplenorm{u - u_h} \lesssim ( \alpha^{\frac12} h^{k_s} + \Delta t^{k_t + \frac12} ) \Vert u \Vert_{H^{k_{\max}+2}(Q)}$.
\end{theorem}
\begin{proof}
 Combine \cref{prop_strang} with $w_h = I_h u^e$ and \cref{lemma_approx}.
\end{proof}
\begin{remark}[Discussion of the error bound]
 In Case D, the bound \cref{eq_apriori_caseD} coincides (for $\alpha = 1$, where $\Delta t^{k_t+1} \leq \Delta t^{k_t + \frac12}$) with the one obtained in \cite{heimann2025discretizationerroranalysishigh} for the continuous-in-space method (with exact geometry handling) and is optimal in the norm $\tripplenorm{\cdot}$, which contains the diffusion-scaled $H^1$-type seminorm $\alpha^{\frac12} |\cdot|_{1,h}$, the $L^2$ norm and the upwind jumps.
 In Case C, the bound \cref{eq_apriori_caseC} coincides with the one obtained in \cite{burman2025transport} for the transport problem ($\alpha = 0$) and is the typical estimate for upwind DG methods for convection dominated problems, cf. \cite{di2011mathematical, ern2021finite3}: it is of order $k+\frac12$ in the $L^2$ norm, which is half an order suboptimal, but optimal for the material derivative error $\Delta t^{\frac12} \Vert \w^\ast \cdo \nabla^\ast (u - u_h) \Vert_Q \lesssim h^{k+\frac12}$, i.e. $\Vert \w^\ast \cdo \nabla^\ast (u - u_h) \Vert_Q \lesssim h^{k}$. The bound is robust in the sense that all constants in Case C are independent of $\alpha \geq 0$; the restriction $\alpha \lesssim \Vert \w \Vert_\infty h$ is only needed to simplify the rate to $h^{k+\frac12}$.
 For $\alpha \gtrsim \Vert \w \Vert_\infty h$ under \cref{eq:ass_caseC}, the assumptions of Case D hold as well, and \cref{eq_apriori_caseD} is the sharper bound, as it does not contain the contribution $h^{k_s + \frac12}$ of the spatial upwind jumps.
\end{remark}
We finally note that by virtue of the problem specific norm the following bounds follow directly from \cref{thm_apriori}.
\begin{corollary} \label{cor_final_bound}
 Under the assumptions of \cref{thm_apriori}, in Case C with $\alpha \lesssim \Vert \w \Vert_\infty h$,
 \begin{align}
  \Vert u - u_h \Vert_Q + \Vert (u - u_h)_- \Vert_{\Omega^N} + h^{\frac12} \Vert \w^\ast \cdo \nabla^\ast (u - u_h) \Vert_Q &\lesssim h^{k + \frac12} \Vert u \Vert_{H^{k_{\max}+2}(Q)},
 \end{align}
 and in Case D
 \begin{align}
  \Vert u - u_h \Vert_Q + \Vert (u - u_h)_- \Vert_{\Omega^N} &+ \alpha^{\frac12} | u - u_h |_{1,h} + \Delta t^{\frac12} \Vert \partial_t (u - u_h) \Vert_Q \nonumber \\ &\lesssim \big( \alpha^{\frac12} ( h^{k_s} + \Delta t^{k_t+1} ) + \Delta t^{k_t + \frac12} \big) \Vert u \Vert_{H^{k_{\max}+2}(Q)}.
 \end{align}
\end{corollary}

\section{Numerical examples} \label{sec:numex}
\input{sec_numex}

\paragraph*{Acknowledgments}
 The authors were supported by the Swedish Research Council under grant no. 2021-06594 while in residence at Institut Mittag-Leffler in Djursholm, Sweden during the fall semester of 2025, whilst conceptualising the method. FH acknowledges support by EPSRC under grant EP/X042650/1.
\paragraph*{Statement on the use of AI}
 The authors conceptualised the method in autumn 2025 and developed a draft analysis without AI assistance. Later, the completion of this proof scaffolding was executed by Anthropic Claude Fable, which also contributed as a coding assistant to the section on numerical examples under close supervision of the authors. Later, OpenAI ChatGPT Astra and Claude Opus were used for reviewing. The authors manually reviewed the analysis with particular attention to the most novel core estimates regarding the upwind in space terms in the inf-sup stability argument, and take full responsibility for all parts of analysis and numerical examples.

 Finally, Claude Fable was used to formalise the analysis of \cref{sect_stab_analysis,sect_a_pr_est} in Lean 4 with Mathlib (folder \texttt{StDgdg} of the repository \url{https://github.com/fabianheimann/st_dgdg_release}). There, the stability, well-posedness and a priori error results are derived from an explicit list of discrete inequalities and identities (\texttt{Setting.lean}), and these derivations are checked by the Lean kernel. The list itself is assumed; within Lean it is verified only for three model problems: an ordinary differential equation, and a stationary and a moving unfitted interval with piecewise linear functions. The authors have checked that the assumptions and main theorems in Lean match those of the paper.

\input{st_article.bbl}
\end{document}

%% file: st_shared_siam.tex
\usepackage{lipsum}
\usepackage{amsfonts}
\usepackage{amsmath}
\usepackage{amssymb}
\usepackage{nicefrac}
\usepackage{graphicx}
\usepackage[caption=false]{subfig}
\usepackage{epstopdf}
\usepackage{algorithmic}
\ifpdf
  \DeclareGraphicsExtensions{.eps,.pdf,.png,.jpg}
\else
  \DeclareGraphicsExtensions{.eps}
\fi

\usepackage[normalem]{ulem}

\newsiamremark{remark}{Remark}
\crefname{remark}{Remark}{Remarks}
\newsiamremark{assumption}{Assumption}
\crefname{assumption}{Assumption}{Assumptions}
\newsiamremark{hypothesis}{Hypothesis}
\crefname{hypothesis}{Hypothesis}{Hypotheses}
\newsiamthm{claim}{Claim}
\newsiamthm{conjecture}{Conjecture}
\crefname{claim}{Claim}{Claims}
\crefname{conjecture}{Conjecture}{Conjectures}
\newsiamthm{attempt}{Attempt}

\usepackage{amsopn}

\usepackage[]{hyperref}
\usepackage[capitalize]{cleveref}

\providecommand{\vertiii}[1]{{\left\vert\kern-0.25ex\left\vert\kern-0.25ex\left\vert #1
\right\vert\kern-0.25ex\right\vert\kern-0.25ex\right\vert}}

\usepackage{tikz}
\usepackage{tikz-3dplot}
\usetikzlibrary{fpu,positioning,shapes,arrows,calc, arrows.meta, patterns, intersections, decorations.pathreplacing,calligraphy}
\usetikzlibrary{shapes,arrows,positioning,fit,calc,spy, colorbrewer}

\usepackage{pgfplots}
\usepgfplotslibrary{colorbrewer,groupplots}

\pgfplotsset{
cycle list/Set1-5,
cycle multiindex* list={
mark list*\nextlist
Set1-5\nextlist
},
}

\usepackage{shellesc}
\usetikzlibrary{external}

\usepackage{booktabs}

\usepackage{graphics}

\usepackage{rotating}
\usepackage{mathtools}

\providecommand{\bi}[1]{\hyperlink{def:bi}{b_{i}}}

\providecommand{\ThnG}[1]{\hyperref[eq:new_st_isoparam_regions]{\mathcal{T}_{h,#1}^{\Gamma}}}

%% file: sec_numex.tex
\newcommand{\numexexport}{\tikzset{external/export next=false}}
\pgfplotsset{cycle list/Set1-6}
\colorlet{numexk6}{Set1-H} %
\pgfplotsset{numex/.style={compat=1.18, thick, mark options={scale=1.2, solid}, legend style={draw=none, fill=none, font=\footnotesize, legend cell align=left}, ymajorgrids, xmajorgrids, grid style={gray!25}, xlabel style={at={(rel axis cs:1,0)}, anchor=north west, xshift=5pt, yshift=-3pt, inner sep=1pt}}}

In this section we illustrate the performance of the unfitted space-time DG method and compare it with its continuous-in-space counterpart. All experiments are implemented in \texttt{ngsxfem} \cite{xfem_joss}, with equal polynomial orders in space and time, $k_s = k_t = k$, and the parameters $\gamma_J = 0.05$ and $\lambda = 10 k^2$. Unlike the analyses of the continuous-in-space methods in \cite{heimann2025discretizationerroranalysishigh, burman2025transport}, the analysis in this paper assumes exact geometry handling. In the computations we realise this approximately by the same construction as in these works, the space-time isoparametric mapping $\Theta_h$ of order $k$ in space and time from \cite{HLP2022, heimann2024geometrically}, for which the geometry error is of the same asymptotic order as the discretisation error; high-order quadrature on implicitly defined domains as in \cite{saye15} would be an alternative. The upwind fluxes then involve the velocity relative to the mesh deformation, $(\w - \partial_t \Theta_h) \cdo \n$. The linear systems on each time slab are solved with the sparse direct solver \texttt{pardiso} (Intel MKL), on the finest level in \cref{ssec:numex_3d} with \texttt{MUMPS} \cite{mumps01} on 16 MPI processes via PETSc \cite{petsc-user-ref}. The scripts and the data underlying the figures and tables of this section are available at \url{https://github.com/fabianheimann/st_dgdg_release}.

As a reference we use the method with the same DG-in-time discretisation, geometry handling and Ghost penalty scaling $c_J = \gamma_J (\alpha h^{-2} + \Vert \w \Vert_\infty h^{-1})$, but continuous finite elements of order $k$ in space and the Ghost penalty acting on all facets between active elements, as in \cite{burman2025transport}, where it also stabilises the convection in the spirit of continuous interior penalty methods. We denote it by CG. For comparison in \cref{ssec:numex_conv}, we also consider the variant CG$_\Gamma$, in which, as for the DG method, the Ghost penalty acts only on the facets of cut elements. It stabilises the cut configurations, but not the convection, cf.\ \cite{burman2025transport}.

\begin{figure}
 \centering
 \tikzset{lab/.style={font=\footnotesize, inner sep=1pt, text=black, fill=white, fill opacity=0.75, text opacity=1, rounded corners=1pt},
          pin/.style={black!70, thin}, panel/.style={font=\footnotesize\bfseries, anchor=north west, inner sep=0pt},
          ptitle/.style={font=\footnotesize, anchor=north east, inner sep=0pt}}
 \numexexport
 \begin{tikzpicture}
  \node[anchor=south west, inner sep=0] (img) at (0,0) {\includegraphics[height=3.0cm]{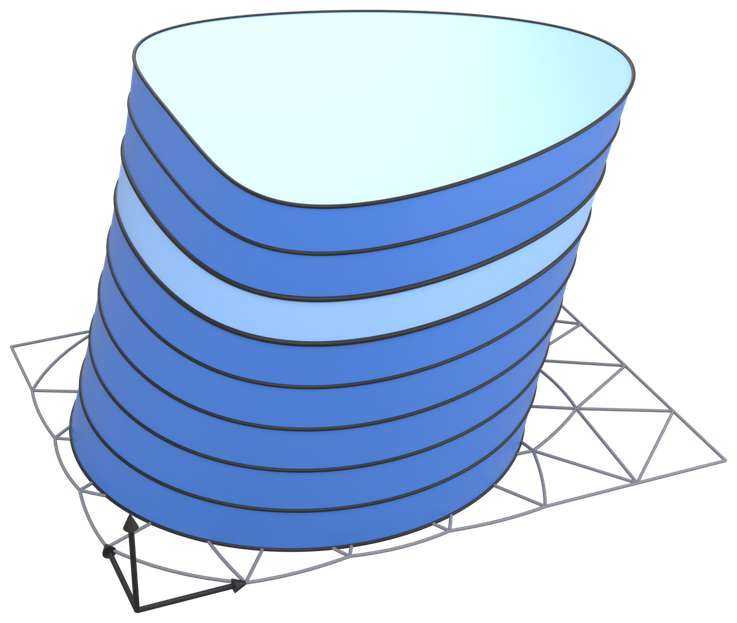}};
  \begin{scope}[x={(img.south east)}, y={(img.north west)}]
   \input{graphics/spacetime-kite-a.anchors.tex}
   \node[panel] at (0,1) {(a)};
   \node[lab, anchor=west] at (a-x) {$x$};
   \node[lab, anchor=south east] at (a-y) {$y$};
   \node[lab, anchor=south] at (a-t) {$t$};
   \node[font=\small, text=white] at (a-Q) {$Q$};
   \draw[pin] (a-In) -- ++(0.08,0.0) node[lab, anchor=west] {$Q^n$};
  \end{scope}
 \end{tikzpicture}\hfill
 \numexexport
 \begin{tikzpicture}
  \node[anchor=south west, inner sep=0] (img) at (0,0) {\includegraphics[height=3.0cm]{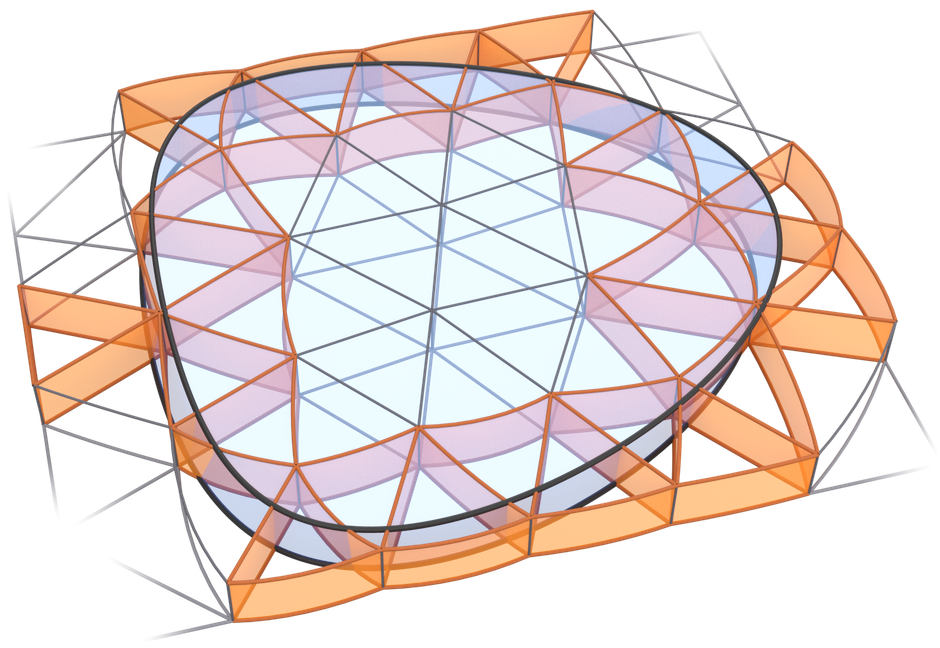}};
  \begin{scope}[x={(img.south east)}, y={(img.north west)}]
   \node[panel] at (0,1) {(b)};
   \node[ptitle] at (1,1) {DG, CG$_\Gamma$};
  \end{scope}
 \end{tikzpicture}\hfill
 \numexexport
 \begin{tikzpicture}
  \node[anchor=south west, inner sep=0] (img) at (0,0) {\includegraphics[height=3.0cm]{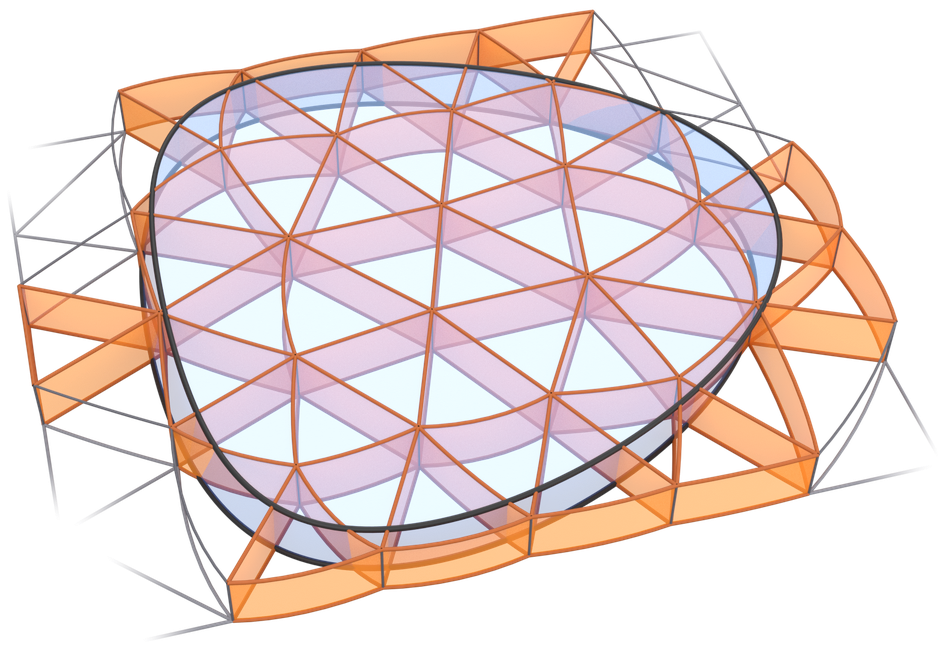}};
  \begin{scope}[x={(img.south east)}, y={(img.north west)}]
   \node[panel] at (0,1) {(c)};
   \node[ptitle] at (1,1) {CG};
  \end{scope}
 \end{tikzpicture}
 \caption{Kite example on the coarsest level $i = 1$ ($h = 0.45$, $8$ time slabs; time axis scaled). (a) Space-time domain $Q$ over the background mesh, with the time levels $t_n$ and the slab $Q^n$, $n = 6$, highlighted. (b), (c) Active prisms of this slab under the isoparametric mapping $\Theta_h$ ($k = 2$, true scale) and the space-time facets $F \times I_n$ with Ghost penalty (orange): $F \in \mathcal{F}^n_{\text{gp}}$ for DG and CG$_\Gamma$, cf.\ \cref{eq:def_gp_facets}, and all facets between active elements for CG.}
 \label{fig:spacetime_kite}
\end{figure}
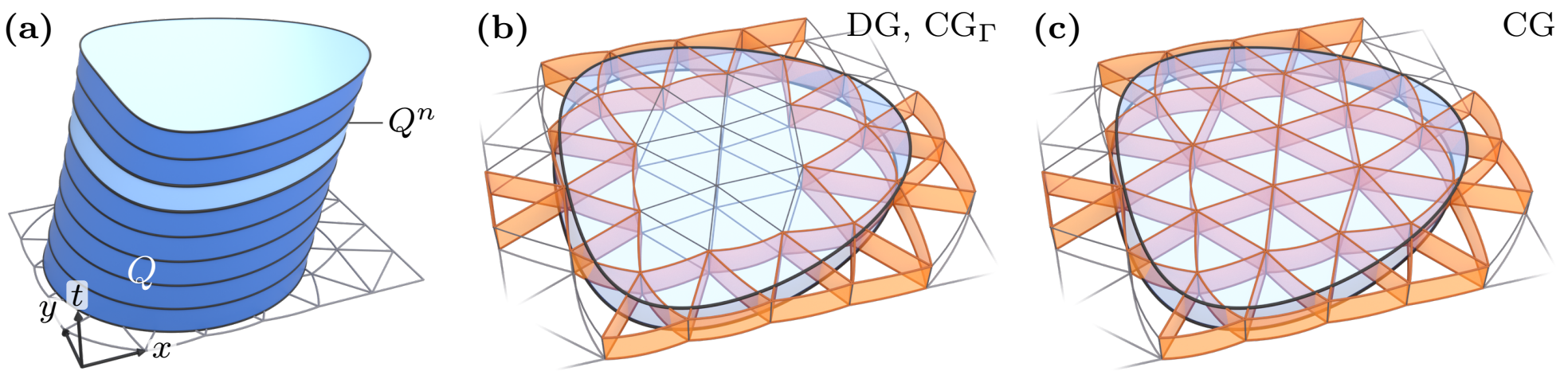

We use the kite geometry from \cite{HLP2022}: a circle that is deformed by a shear flow into a kite-shaped domain,
\begin{equation*}
 \phi(x,y,t) = \sqrt{(x - \rho)^2 + y^2} - 1, \qquad \rho(y,t) = (1 - y^2)\, t, \qquad \w = (1 - y^2, 0),
\end{equation*}
$t \in [0, 0.5]$, in the background domain $\tilde \Omega = [-1.05, 1.55] \times [-1.05, 1.05]$. As $\w$ transports $\phi$, the domain moves with $\w$ as assumed in \cref{sec:model_problem}. On refinement level $i$ we use unstructured simplicial meshes with $h = 0.9 \cdot 2^{-i}$ and $\Delta t = 0.125 \cdot 2^{-i}$, i.e.\ space and time are refined simultaneously and \cref{eq:ass_caseC} holds. The data $f$, $u_0$ are chosen such that $u = \cos(\pi r) \sin(\pi t)$, $r = \sqrt{(x - \rho)^2 + y^2}$, is the exact solution. \Cref{fig:spacetime_kite} illustrates the space-time domain, the active space-time mesh and the facets with Ghost penalty on the coarsest level. We report the errors $\Vert u - u_h \Vert_{\Omega(T)} := \Vert u(T) - (u_h)_-(T) \Vert_{L^2(\Omega(T))}$ and the broken $H^1$ seminorm error $|u - u_h|_{1,h, Q}$ (volume contributions only).

\subsection{Convergence} \label{ssec:numex_conv}
\begin{figure}
 \centering
 {\footnotesize \hypersetup{linkcolor=black}\numexexport\ref{numexlegconv}} \\[1ex]
 \numexexport
 \begin{tikzpicture}
 \begin{semilogyaxis}[numex, width=0.45\textwidth, height=0.42\textwidth, xlabel={$i$}, ylabel={$\Vert u - u_h \Vert_{\Omega(T)}$}, legend to name=numexlegconv, legend columns=6, legend style={/tikz/every even column/.append style={column sep=6pt}}, xtick={1,...,7}, ymin=1e-10]
 \addplot+[numexk6, mark=o, forget plot] table[x=i,y=l2final] {num_exp/out/conv_kite_DG_ks6_kt6_a1_g0.05_c1.dat}; %
 \addplot+[Set1-A, mark=*] table[x=i,y=l2final] {num_exp/out/conv_kite_DG_ks1_kt1_a1_g0.05_c1.dat}; \addlegendentry{$k=1$}
 \addplot+[Set1-B, mark=square*] table[x=i,y=l2final] {num_exp/out/conv_kite_DG_ks2_kt2_a1_g0.05_c1.dat}; \addlegendentry{$k=2$}
 \addplot+[Set1-C, mark=triangle*] table[x=i,y=l2final] {num_exp/out/conv_kite_DG_ks3_kt3_a1_g0.05_c1.dat}; \addlegendentry{$k=3$}
 \addplot+[Set1-D, mark=diamond*] table[x=i,y=l2final] {num_exp/out/conv_kite_DG_ks4_kt4_a1_g0.05_c1.dat}; \addlegendentry{$k=4$}
 \addplot+[Set1-E, mark=pentagon*] table[x=i,y=l2final] {num_exp/out/conv_kite_DG_ks5_kt5_a1_g0.05_c1.dat}; \addlegendentry{$k=5$}
 \addlegendimage{numexk6, mark=o} \addlegendentry{$k=6$}
 \addplot[gray, dotted, domain=1:7, forget plot] {0.05*2^(-2*x)};
 \addplot[gray, dotted, domain=1:6, forget plot] {0.02*2^(-3*x)};
 \addplot[gray, dotted, domain=1:5, forget plot] {0.03*2^(-4*x)};
 \addplot[gray, dotted, domain=1:4, forget plot] {0.01*2^(-5*x)};
 \addplot[gray, dotted, domain=1:4, forget plot] {0.006*2^(-6*x)};
 \addplot[gray, dotted, domain=1:3, forget plot] {0.003*2^(-7*x)};
 \end{semilogyaxis}
 \end{tikzpicture}
 \numexexport
 \begin{tikzpicture}
 \begin{semilogyaxis}[numex, width=0.45\textwidth, height=0.42\textwidth, xlabel={$i$}, ylabel={$|u - u_h|_{1,h,Q}$}, xtick={1,...,7}, ymin=1e-9]
 \addplot+[numexk6, mark=o, forget plot] table[x=i,y=h1] {num_exp/out/conv_kite_DG_ks6_kt6_a1_g0.05_c1.dat}; %
 \addplot+[Set1-A, mark=*] table[x=i,y=h1] {num_exp/out/conv_kite_DG_ks1_kt1_a1_g0.05_c1.dat}; \addlegendentry{$k=1$}
 \addplot+[Set1-B, mark=square*] table[x=i,y=h1] {num_exp/out/conv_kite_DG_ks2_kt2_a1_g0.05_c1.dat}; \addlegendentry{$k=2$}
 \addplot+[Set1-C, mark=triangle*] table[x=i,y=h1] {num_exp/out/conv_kite_DG_ks3_kt3_a1_g0.05_c1.dat}; \addlegendentry{$k=3$}
 \addplot+[Set1-D, mark=diamond*] table[x=i,y=h1] {num_exp/out/conv_kite_DG_ks4_kt4_a1_g0.05_c1.dat}; \addlegendentry{$k=4$}
 \addplot+[Set1-E, mark=pentagon*] table[x=i,y=h1] {num_exp/out/conv_kite_DG_ks5_kt5_a1_g0.05_c1.dat}; \addlegendentry{$k=5$}
 \legend{}
 \addplot[gray, dotted, domain=1:7, forget plot] {0.8*2^(-x)};
 \addplot[gray, dotted, domain=1:6, forget plot] {0.5*2^(-2*x)};
 \addplot[gray, dotted, domain=1:5, forget plot] {0.2*2^(-3*x)};
 \addplot[gray, dotted, domain=1:4, forget plot] {0.1*2^(-4*x)};
 \addplot[gray, dotted, domain=1:4, forget plot] {0.05*2^(-5*x)};
 \addplot[gray, dotted, domain=1:3, forget plot] {0.03*2^(-6*x)};
 \end{semilogyaxis}
 \end{tikzpicture}
 \\[1ex]
 \numexexport
 \begin{tikzpicture}
 \begin{semilogyaxis}[numex, width=0.45\textwidth, height=0.42\textwidth, xlabel={$i$}, ylabel={$\Vert u - u_h \Vert_{\Omega(T)}$}, xtick={1,...,6}, xmin=0.7, xmax=6.3, ymin=1e-10]
 \addplot+[numexk6, mark=o, forget plot] table[x=i,y=l2final] {num_exp/out/conv_kite_DG_ks6_kt6_a0_g0.05_c1.dat};
 \addplot+[Set1-A, mark=*, forget plot] table[x=i,y=l2final] {num_exp/out/conv_kite_DG_ks1_kt1_a0_g0.05_c1.dat};
 \addplot+[Set1-B, mark=square*, forget plot] table[x=i,y=l2final] {num_exp/out/conv_kite_DG_ks2_kt2_a0_g0.05_c1.dat};
 \addplot+[Set1-C, mark=triangle*, forget plot] table[x=i,y=l2final] {num_exp/out/conv_kite_DG_ks3_kt3_a0_g0.05_c1.dat};
 \addplot+[Set1-D, mark=diamond*, forget plot] table[x=i,y=l2final] {num_exp/out/conv_kite_DG_ks4_kt4_a0_g0.05_c1.dat};
 \addplot+[Set1-E, mark=pentagon*, forget plot] table[x=i,y=l2final] {num_exp/out/conv_kite_DG_ks5_kt5_a0_g0.05_c1.dat};
 \addplot[gray, dotted, domain=1:6] {0.15*2^(-2*x)};
 \addplot[gray, dotted, domain=1:4] {0.006*2^(-6*x)};
 \addplot[gray, dotted, domain=1:3] {0.003*2^(-7*x)};
 \addplot[gray, dotted, domain=1:5, forget plot] {0.02*2^(-3*x)};
 \addplot[gray, dotted, domain=1:5, forget plot] {0.03*2^(-4*x)};
 \addplot[gray, dotted, domain=1:4, forget plot] {0.01*2^(-5*x)};
 \end{semilogyaxis}
 \end{tikzpicture}
 \numexexport
 \begin{tikzpicture}
 \begin{semilogyaxis}[numex, width=0.45\textwidth, height=0.42\textwidth, xlabel={$i$}, ylabel={$\Vert \w^\ast \cdo \nabla^\ast (u - u_h) \Vert_{Q}$}, xtick={1,...,6}, xmin=0.7, xmax=6.3, ymin=1e-9]
 \addplot+[numexk6, mark=o, forget plot] table[x=i,y=mat] {num_exp/out/conv_kite_DG_ks6_kt6_a0_g0.05_c1.dat};
 \addplot+[Set1-A, mark=*, forget plot] table[x=i,y=mat] {num_exp/out/conv_kite_DG_ks1_kt1_a0_g0.05_c1.dat};
 \addplot+[Set1-B, mark=square*, forget plot] table[x=i,y=mat] {num_exp/out/conv_kite_DG_ks2_kt2_a0_g0.05_c1.dat};
 \addplot+[Set1-C, mark=triangle*, forget plot] table[x=i,y=mat] {num_exp/out/conv_kite_DG_ks3_kt3_a0_g0.05_c1.dat};
 \addplot+[Set1-D, mark=diamond*, forget plot] table[x=i,y=mat] {num_exp/out/conv_kite_DG_ks4_kt4_a0_g0.05_c1.dat};
 \addplot+[Set1-E, mark=pentagon*, forget plot] table[x=i,y=mat] {num_exp/out/conv_kite_DG_ks5_kt5_a0_g0.05_c1.dat};
 \addplot[gray, dotted, domain=1:6] {0.5*2^(-x)};
 \addplot[gray, dotted, domain=1:4] {0.006*2^(-5*x)};
 \addplot[gray, dotted, domain=1:3] {0.004*2^(-6*x)};
 \addplot[gray, dotted, domain=1:5, forget plot] {0.2*2^(-2*x)};
 \addplot[gray, dotted, domain=1:5, forget plot] {0.05*2^(-3*x)};
 \addplot[gray, dotted, domain=1:4, forget plot] {0.015*2^(-4*x)};
 \end{semilogyaxis}
 \end{tikzpicture}
 \caption{DG method under simultaneous refinement in space and time. Top: $\alpha = 1$, $k = 1, \dots, 6$, $L^2$ error at final time (left, dotted: $\mathcal{O}(h^{k+1})$) and broken $H^1$ seminorm error (right, dotted: $\mathcal{O}(h^{k})$). Bottom: $\alpha = 0$, $k = 1, \dots, 6$, $L^2$ error at final time (left, dotted: $\mathcal{O}(h^{k+1})$) and material derivative error (right, dotted: $\mathcal{O}(h^{k})$).}
 \label{fig:conv}
\end{figure}
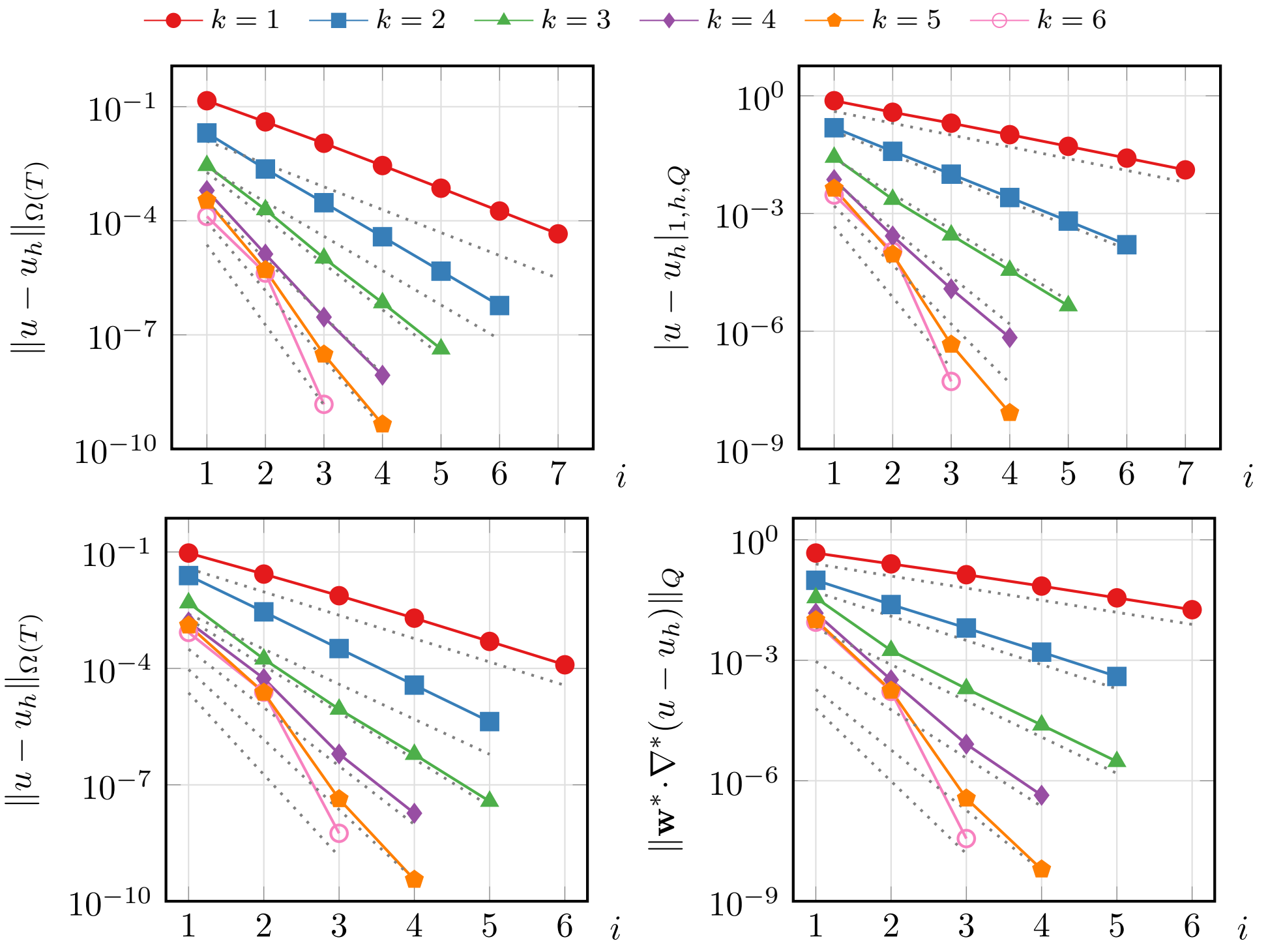

\Cref{fig:conv} shows the errors of the DG method under refinement. For $\alpha = 1$ and $k = 1, \dots, 6$, we observe the optimal orders $\mathcal{O}(h^{k+1})$ in the $L^2$ norm and $\mathcal{O}(h^k)$ in the broken $H^1$ seminorm, the latter as predicted by \cref{thm_apriori} (Case D). For $\alpha = 0$, the material derivative error converges with order $k$, as predicted by \cref{thm_apriori} (Case C), and the $L^2$ error with order $k+1$, half an order better than the bound. For $k \geq 5$ the coarsest levels are not yet in the asymptotic regime, so that the observed orders on the available levels are higher.

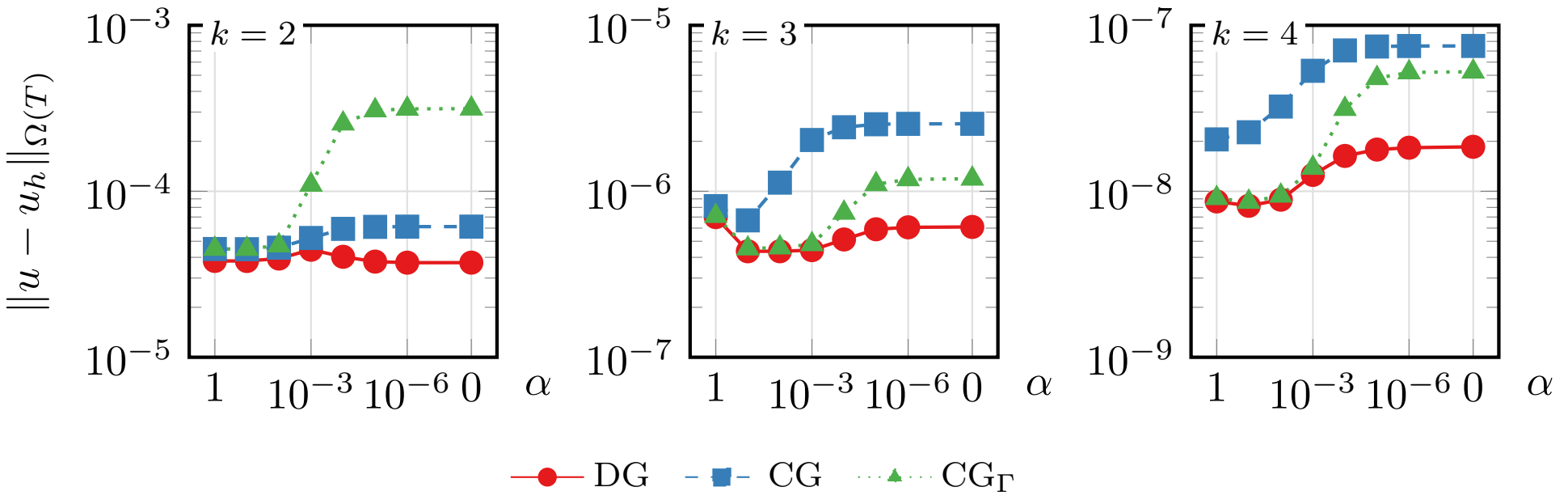
\begin{figure}
 \centering

 \numexexport
 \begin{tikzpicture}
 \begin{loglogaxis}[numex, width=0.315\textwidth, height=0.33\textwidth, xlabel={$\alpha$}, title={$k = 2$}, title style={at={(0.04,0.96)}, anchor=north west, font=\footnotesize, fill=white, inner sep=1.5pt}, x dir=reverse, xtick={1,1e-3,1e-6,1e-8}, xticklabels={$1$,$10^{-3}$,$10^{-6}$,$0$}, xminorticks=false, ymin=1e-5, ymax=1e-3, ylabel={$\Vert u - u_h \Vert_{\Omega(T)}$}, legend to name=numexlegpeclet, legend columns=3, legend style={/tikz/every even column/.append style={column sep=6pt}}]
 \IfFileExists{num_exp/out/peclet_k2_i4.dat}{
 \addplot[Set1-A, mark=*] table[x=alpha, y=DG] {num_exp/out/peclet_k2_i4.dat}; \addlegendentry{DG}
 \addplot[Set1-B, mark=square*, dashed] table[x=alpha, y=CGall] {num_exp/out/peclet_k2_i4.dat}; \addlegendentry{CG}
 \addplot[Set1-C, mark=triangle*, dotted] table[x=alpha, y=CG] {num_exp/out/peclet_k2_i4.dat}; \addlegendentry{CG$_\Gamma$}}{}
 \end{loglogaxis}
 \end{tikzpicture}
 \numexexport
 \begin{tikzpicture}
 \begin{loglogaxis}[numex, width=0.315\textwidth, height=0.33\textwidth, xlabel={$\alpha$}, title={$k = 3$}, title style={at={(0.04,0.96)}, anchor=north west, font=\footnotesize, fill=white, inner sep=1.5pt}, x dir=reverse, xtick={1,1e-3,1e-6,1e-8}, xticklabels={$1$,$10^{-3}$,$10^{-6}$,$0$}, xminorticks=false, ymin=1e-7, ymax=1e-5, legend to name=numexlegpecletunusedb]
 \IfFileExists{num_exp/out/peclet_k3_i4.dat}{
 \addplot[Set1-A, mark=*] table[x=alpha, y=DG] {num_exp/out/peclet_k3_i4.dat}; \addlegendentry{DG}
 \addplot[Set1-B, mark=square*, dashed] table[x=alpha, y=CGall] {num_exp/out/peclet_k3_i4.dat}; \addlegendentry{CG}
 \addplot[Set1-C, mark=triangle*, dotted] table[x=alpha, y=CG] {num_exp/out/peclet_k3_i4.dat}; \addlegendentry{CG$_\Gamma$}}{}
 \end{loglogaxis}
 \end{tikzpicture}
 \numexexport
 \begin{tikzpicture}
 \begin{loglogaxis}[numex, width=0.315\textwidth, height=0.33\textwidth, xlabel={$\alpha$}, title={$k = 4$}, title style={at={(0.04,0.96)}, anchor=north west, font=\footnotesize, fill=white, inner sep=1.5pt}, x dir=reverse, xtick={1,1e-3,1e-6,1e-8}, xticklabels={$1$,$10^{-3}$,$10^{-6}$,$0$}, xminorticks=false, ymin=1e-9, ymax=1e-7, legend to name=numexlegpecletunusedc]
 \IfFileExists{num_exp/out/peclet_k4_i4.dat}{
 \addplot[Set1-A, mark=*] table[x=alpha, y=DG] {num_exp/out/peclet_k4_i4.dat}; \addlegendentry{DG}
 \addplot[Set1-B, mark=square*, dashed] table[x=alpha, y=CGall] {num_exp/out/peclet_k4_i4.dat}; \addlegendentry{CG}
 \addplot[Set1-C, mark=triangle*, dotted] table[x=alpha, y=CG] {num_exp/out/peclet_k4_i4.dat}; \addlegendentry{CG$_\Gamma$}}{}
 \end{loglogaxis}
 \end{tikzpicture}
 \\[1ex] {\footnotesize \hypersetup{linkcolor=black}\numexexport\ref{numexlegpeclet}}
 \caption{Dependence on the diffusion coefficient on the fixed level $i = 4$: $L^2$ error at final time of the DG method, the CG method (Ghost penalty on all facets) and the variant CG$_\Gamma$ (Ghost penalty on the facets of cut elements only) for $\alpha \in \{1, 10^{-1}, \dots, 10^{-6}, 0\}$; all panels span two decades. For $k = 2$ and $\alpha = 0$, the observed order of CG$_\Gamma$ under refinement drops to two.}
 \label{fig:peclet}
\end{figure}

\paragraph{Dependence on the diffusion coefficient}
\Cref{fig:peclet} shows the errors on the fixed level $i = 4$ for $\alpha$ from $1$ down to $0$ and $k \in \{2, 3, 4\}$. The errors of the DG method change by at most a factor of about two over the whole range of $\alpha$. The errors of the variant CG$_\Gamma$, which does not stabilise the convection, increase once convection dominates, whereas the CG method remains robust for $\alpha \to 0$, but is for $k \in \{3, 4\}$ less accurate than the DG method by a factor of up to four. Hence, the DG method is uniformly accurate over the whole range of $\alpha$ with one fixed set of parameters.

\paragraph{Stabilisation parameters}
As known for the continuous-in-space method \cite{HLP2022}, the Ghost penalty parameter can be chosen from a wide range: for $k = 4$ and $\gamma_J \in [5 \cdot 10^{-4}, 5]$ we observe optimal convergence for both $\alpha = 1$ and $\alpha = 0$, with errors on the finest level within a factor of about four (not shown). For the interior penalty, we found $\lambda \geq k^2$ to be sufficient, with practically identical errors up to $\lambda = 250 k^2$ (not shown).

\paragraph{Comparison with the continuous-in-space method}
On the same mesh, the DG method has more active degrees of freedom than the CG method, asymptotically by the factor $(k+1)(k+2) / k^2$ for simplicial meshes in 2D. We observe these ratios (not shown), from about $6$ for $k = 1$ down to $1.5$ for $k = 6$, and practically the same ones for the numbers of non-zero entries of the system matrices. Hence, the relative overhead of the DG method decreases for higher polynomial degree.

\subsection{Three space dimensions} \label{ssec:numex_3d}
The method and its implementation carry over to three space dimensions without changes. To demonstrate this, \cref{tab:3d} shows the errors of the DG method for a moving sphere: $\phi(x,t) = \Vert x - \rho(t) e_2 \Vert_2 - 0.5$, $\rho(t) = \frac{1}{\pi} \sin(2 \pi t)$, $\w = \rho'(t) e_2$, $t \in [0, 0.5]$, in the background domain $\tilde \Omega = [-0.6, 0.6] \times [-1, 1] \times [-0.6, 0.6]$, with $h = 0.5 \cdot 2^{-i}$, $\Delta t = 0.125 \cdot 2^{-i}$ and the exact solution $u = \cos(2 \pi r) \sin(\pi t)$, $r = \Vert x - \rho(t) e_2 \Vert_2$. For $k = 3$, the observed orders on the coarser levels exceed $k + 1$; the best approximation of $u(T)$ on these meshes shows the same pre-asymptotic behaviour, and on the finest level the orders drop to values between $k + \frac12$ and $k + 1$. For corresponding computations with continuous-in-space methods we refer to \cite{HLP2022} ($\alpha = 1$) and \cite{burman2025transport} ($\alpha = 0$).

\begin{table}
 \centering \footnotesize
 \caption{DG method in three space dimensions (moving sphere): active degrees of freedom per time slab, $L^2$ error at final time and observed order (eoc) for $\alpha = 1$ and $\alpha = 0$.}
 \label{tab:3d}
 \begin{tabular}{ccr cc cc}
 \toprule
 & & & \multicolumn{2}{c}{$\alpha = 1$} & \multicolumn{2}{c}{$\alpha = 0$} \\
 \cmidrule(lr){4-5} \cmidrule(lr){6-7}
 $k$ & $i$ & dofs & $\Vert u - u_h \Vert_{\Omega(T)}$ & eoc & $\Vert u - u_h \Vert_{\Omega(T)}$ & eoc \\
 \midrule
 2 & 0 & 1770 & $7.44\cdot 10^{-1}$ & -- & $1.06\cdot 10^{-1}$ & -- \\
  & 1 & 14220 & $5.57\cdot 10^{-2}$ & 3.7 & $1.16\cdot 10^{-2}$ & 3.2 \\
  & 2 & 82680 & $3.99\cdot 10^{-3}$ & 3.8 & $1.63\cdot 10^{-3}$ & 2.8 \\
  & 3 & 506940 & $3.02\cdot 10^{-4}$ & 3.7 & $2.27\cdot 10^{-4}$ & 2.8 \\
 \midrule
 3 & 0 & 4720 & $7.07\cdot 10^{-1}$ & -- & $1.44\cdot 10^{-1}$ & -- \\
  & 1 & 37920 & $2.70\cdot 10^{-2}$ & 4.7 & $4.18\cdot 10^{-3}$ & 5.1 \\
  & 2 & 220480 & $1.11\cdot 10^{-3}$ & 4.6 & $1.12\cdot 10^{-4}$ & 5.2 \\
  & 3 & 1351840 & $9.98\cdot 10^{-5}$ & 3.5 & $8.50\cdot 10^{-6}$ & 3.7 \\
 \bottomrule
 \end{tabular}
\end{table}

\subsection{Transport of a steep layer} \label{ssec:numex_ring}
Upwind DG methods are known to perform particularly well for convection dominated problems with steep, barely resolved layers, where they keep spurious oscillations small and local to the layer, cf.\ e.g.\ \cite{di2011mathematical}. To examine this in the unfitted space-time setting, we keep the kite geometry and the velocity field $\w$ from \cref{ssec:numex_conv} and choose $\alpha = 0$, $f = 0$ and $u_0 = g(r)$, $g(r) = (1 + \exp(2 (r - 0.5) / \delta))^{-1}$, so that the exact solution $u = g(r)$ is a ring-shaped layer of width $\sim \delta$ which is transported by $\w$.

\begin{figure}
 \centering \numexexport
 \begin{tikzpicture}
 \begin{axis}[numex, width=0.42\textwidth, height=0.32\textwidth, xlabel={$x$}, ylabel={$u_h(x, 0, T)$}, xmin=-0.4, xmax=1.4, ymin=-0.0002, ymax=0.0002, ytick={-0.0002,-0.0001,0,0.0001,0.0002}, yticklabels={$-2 \cdot 10^{-4}$,$-10^{-4}$,$0$,$10^{-4}$,$2 \cdot 10^{-4}$}, scaled y ticks=false, legend style={at={(0.5,0.5)}, anchor=center}, restrict y to domain*=-0.0004:0.0004]
 \IfFileExists{num_exp/out/ring_profile_DG.dat}{\addplot[black, thin] table[x=x, y=u] {num_exp/out/ring_profile_DG.dat}; \addlegendentry{exact}
 \addplot[Set1-A, thick] table[x=x, y=uh] {num_exp/out/ring_profile_DG.dat}; \addlegendentry{DG}}{}
 \IfFileExists{num_exp/out/ring_profile_CGall.dat}{\addplot[Set1-B, thick, dashed] table[x=x, y=uh] {num_exp/out/ring_profile_CGall.dat}; \addlegendentry{CG, $\gamma_J = 0.05$}}{}
 \IfFileExists{num_exp/out/ring_profile_CGall_g0.005.dat}{\addplot[Set1-C, thin, densely dotted] table[x=x, y=uh] {num_exp/out/ring_profile_CGall_g0.005.dat}; \addlegendentry{CG, $\gamma_J = 0.005$}}{}
 \end{axis}
 \end{tikzpicture}
 \numexexport
 \begin{tikzpicture}
 \begin{semilogyaxis}[numex, width=0.42\textwidth, height=0.32\textwidth, xlabel={$x$}, ylabel={$\vert u - u_h \vert (x, 0, T)$}, xmin=-0.4, xmax=1.4, ymin=1e-12, ymax=1]
 \IfFileExists{num_exp/out/ring_profile_CGall_g0.005.dat}{\addplot[Set1-C, thin, densely dotted] table[x=x, y=err] {num_exp/out/ring_profile_CGall_g0.005.dat};}{}
 \IfFileExists{num_exp/out/ring_profile_CGall.dat}{\addplot[Set1-B, thin, dashed] table[x=x, y=err] {num_exp/out/ring_profile_CGall.dat};}{}
 \IfFileExists{num_exp/out/ring_profile_DG.dat}{\addplot[Set1-A, thin] table[x=x, y=err] {num_exp/out/ring_profile_DG.dat};}{}
 \end{semilogyaxis}
 \end{tikzpicture}
 \caption{Steep layer, $\delta = 0.03$, $\alpha = 0$, $k = 4$, along the line $y = 0$ at final time, for equal numbers of degrees of freedom (DG: level $i = 4$; CG: $h \approx 0.041$, $87$ time steps): discrete solutions, zoomed to the plateau $u \approx 0$ (left; the layers are at $x \approx 0$ and $x \approx 1$, where $u$ rises to $1$), and pointwise errors (right, same legend). CG is shown for $\gamma_J = 0.05$ (dashed) and $\gamma_J = 0.005$ (dotted).}
 \label{fig:ring_profile}
\end{figure}
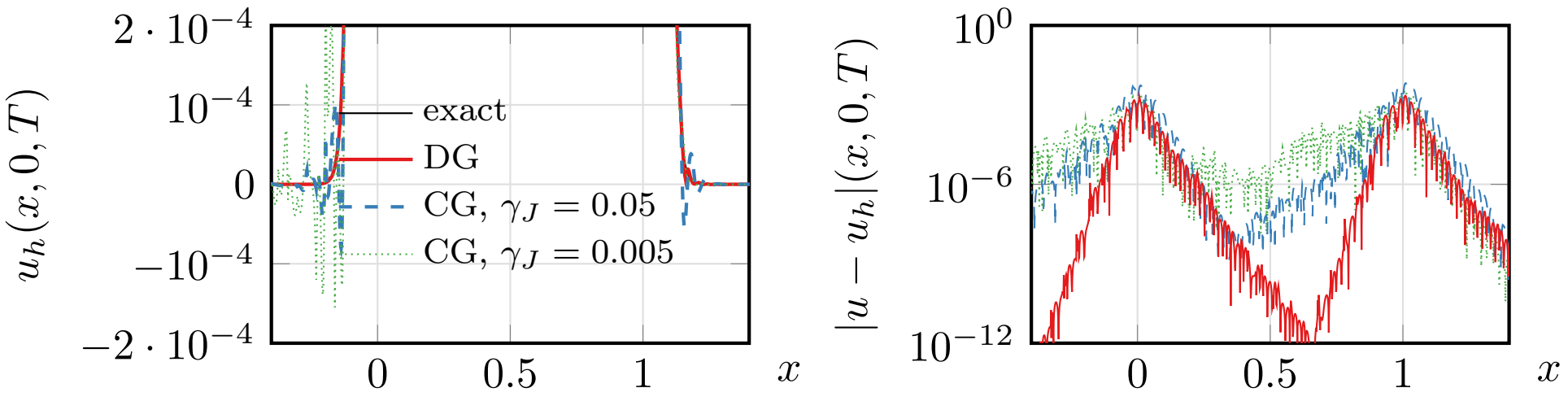

\Cref{fig:ring_profile} compares the DG and the CG method for $\delta = 0.03$ and $k = 4$ at equal numbers of degrees of freedom: the DG method is used on level $i = 4$, where $h \approx 0.056$ is comparable to the layer width, and the CG method on a finer mesh, $h \approx 0.041$ with $87$ time steps, so that these agree up to $1\%$. With $\gamma_J = 0.05$ as before, the CG solution shows over- and undershoots of size $10^{-4}$ next to the layer, about $50$ times larger than those of the DG solution, and its $L^2$ error at final time is larger by a factor of about three. A smaller Ghost penalty parameter, $\gamma_J = 0.005$, brings the $L^2$ error of the CG method close to that of the DG method, but increases its over- and undershoots and its error away from the layer by a further factor of $7$ to $10$. A larger parameter, $\gamma_J = 0.5$, in turn smears the layer and yields a larger $L^2$ error, by a factor of about three compared to $\gamma_J = 0.05$. Hence, for the CG method the choice of $\gamma_J$ is a compromise between the accuracy in the layer and the suppression of spurious oscillations, whereas the DG method yields a small and localised error without such a trade-off.

%% file: graphics/spacetime-kite-a.anchors.tex
\coordinate (a-Q) at (0.3273,0.2803);
\coordinate (a-In) at (0.8086,0.6939);
\coordinate (a-T) at (0.1851,0.8927);
\coordinate (a-x) at (0.3388,0.0663);
\coordinate (a-y) at (0.1423,0.1283);
\coordinate (a-t) at (0.1762,0.1766);